\documentclass{article}

\usepackage{graphicx} 
\usepackage{url}

\usepackage[authoryear]{natbib}
\usepackage[]{hyperref}  %<----modified by Ivan

\usepackage{amsmath}
\usepackage{amssymb}
\usepackage{amsfonts}
\usepackage{amsthm}

\usepackage{multirow}
\usepackage{cancel}

\usepackage[table]{xcolor}
\usepackage{ulem}

\newtheorem{proposition}{Proposition}
\newtheorem{theorem}{Theorem}

\newtheorem{assumption}{Assumption}
\newtheorem{remark}{Remark}

\renewcommand{\tilde}{\widetilde}
\renewcommand{\hat}{\widehat}
\newcommand{\intd}{\,\mathrm{d}}

\usepackage{comment}

\begin{document}

\title{On robustness of spectral R\'{e}nyi divergence}

\author{Tetsuya Takabatake\thanks{
Graduate School of Engineering Science, The University of Osaka, 
           1-3 Machikaneyama, Toyonaka, Osaka 560-8531, Japan. 
           \texttt{t.takabatake.es@osaka-u.ac.jp} 
}
\and
Keisuke Yano\thanks{
The Institute of Statistical Mathematics,
           10-3 Midori-cho, Tachikawa, Tokyo 190‑8562, Japan. 
           \texttt{yano@ism.ac.jp}
}
}

\date{}

\maketitle

\noindent\textbf{Abstract}: This paper studies a specific class of statistical divergences for spectral densities of time series: the spectral $\alpha$-R\'{e}nyi divergences, which include the Itakura--Saito divergence as a limiting case.  
The aim of this paper is to highlight both information-theoretic and statistical properties of spectral $\alpha$-R\'{e}nyi divergences. We reveal 
the connection between the spectral $\alpha$-R\'{e}nyi divergence and the $\gamma$-divergence in robust statistics, and
a variational representation of the spectral $\alpha$-R\'{e}nyi divergence. Inspired by these results suggesting ``robustness" of spectral $\alpha$-R\'{e}nyi divergence, we show that the minimum spectral R\'{e}nyi divergence estimate has a stable optimization path with respect to outliers in the frequency domain, unlike the minimum Itakura–Saito divergence estimator, and thus it delivers more stable estimates, reducing the need for intricate pre-processing.\\


\section{Introduction}

Frequency-domain analysis of time series data has been conducted in many applied fields.
Central to this analysis is the (power) spectral density, a crucial element whose estimation has garnered significant attention. For instance, in seismology, earthquake source parameters are delineated through the spectral parameter estimation \cite{CalderoniandAbercrombie2023,Yoshimitsuetal2023}. 
In geodesy, 
temporal correlation of noises in the Global Navigation Satellite System is captured by spectral parameter estimation
\cite{Langbein_2004}.
In audio signal processing, various methodologies have been developed via spectral density estimation to achieve signal separation \cite{Martin2005}.

Traditional spectral estimation often hinges on the Whittle likelihood maximization \cite{Whittle1953}, which is equivalent to minimizing the Itakura--Saito divergence \cite{ItakuraSaito1968} between a periodogram and a selected class of spectral densities. Broadening the perspective, one can approach such estimations through the lens of spectral divergences, statistical metrics for dissimilarities between two spectral densities. Estimation based on general spectral divergences is comprehensively explored by \cite{Taniguchi1987}.

In this paper, we take another look at a particular class of spectral divergences: the spectral $\alpha$-R\'{e}nyi divergences \cite{Vajda1989,ZhangTaniguchi1995,Kakizawaetal1998,Giletal2013,Griveletal2021}. Notably, 
a limit of the spectral R\'{e}nyi divergence induces the Itakura--Saito divergence. 
In the probabilistic and information-theoretic literature,
this class has been mentioned in the extended discussion of ``Pinsker's information-theoretic justification of the Itakura-Saito distortion measure" \cite{Parzen1993}.
In the statistical literature,
the robustness of the spectral R\'{e}nyi divergence against spectral peaks is pointed out \cite{ZhangTaniguchi1995},
which has been utilized in clustering analysis of time series \cite{Kakizawaetal1998,Hirukawa2005}.

The primary aim of this paper is to elucidate new properties of this class that shed light on the robustness in the spectral parameter estimation.
Specifically, we show that (i) the $\gamma$-divergence \cite{FujisawaEguchi2008} between probability densities of the processes leads to the spectral R\'{e}nyi divergence,
and
(ii) the spectral R\'{e}nyi divergence has a variational representation.
Inspired by these new results, we explore (iii) an additional robustness property of the minimum spectral R\'{e}nyi divergence estimate.
The minimum spectral R\'{e}nyi divergence estimate has been shown to be robust against time series outliers in the frequency domain \cite{ZhangTaniguchi1995,Kakizawaetal1998,Hirukawa2005}.
Expanding on these findings, 
we show that the optimization path in obtaining the minimum spectral R\'{e}nyi divergence estimate
is stable in the presence of time series outliers, a quality that the minimum Itakura--Saito divergence estimator does not have.

Let us explain the implications of our results.
The first result bridges spectral and probabilistic divergences, providing an information-theoretic justification of the spectral R\'{e}nyi divergence.
The second result gives a decomposition of spectral R\'{e}nyi divergence in terms of the Itakura--Saito divergence. The first two results contribute to understanding the robustness of the spectral R\'{e}nyi divergence.
The third result has both theoretical and practical implications.
Theoretically, it gives a novel connection between robust spectral analysis (c.f., \cite{Maronna2019robust}) and optimization theory.
Developing robust statistics from an algorithmic or computational perspective has been a recent topic of discussion (c.f., \cite{Diakonikolas_Kane_2023}).
Our results provide an additional perspective on robust statistics from a computational viewpoint, focusing on the stability of optimization paths.
Practically, 
outliers in the frequency domain often emerge due to insufficient detrending \cite{HeydeandDai1996,Iacone2010,McCloskeyandPerron2013}. Yet, when using the Itakura–Saito divergence, practitioners should exercise greater caution regarding detrending and other pre-processing steps. In comparison, spectral R\'{e}nyi divergences offer more stable estimation results  without necessitating elaborate pre-processing.
This aspect is further supported by an applied study \cite{kano2025spatiotemporal}.

The structure of this paper is as follows.
Section \ref{sec:renyi} delivers new formulae related to the spectral R\'{e}nyi divergence 
(Theorems \ref{thm: Gamma meets Renyi} and \ref{thm: variational rep e}).
Section \ref{sec:outlier} discusses the stability of optimization paths of the spectral R\'{e}nyi divergence minimization (Theorem \ref{thm: success Renyi}). Section \ref{sec:numerical} presents thorough numerical studies. Proofs of the results are presented in Appendices.

\section{Spectral R\'{e}nyi divergence}\label{sec:renyi}

This section introduces a class of spectral R\'{e}nyi divergences and their properties.
Let $\mathcal{S}:=\{S:[-\pi,\pi]\to[0,\infty] \mid \int_{-\pi}^{\pi} S(\omega)\intd \omega<\infty, \, S(\omega)>0 \,\mbox{and}\, S(\omega)=S(-\omega) \, \text{for all $\omega\in (0,\pi)$}\}$.
Fix $S\in\mathcal{S}$ and $\tilde{S}\in\mathcal{S}$ arbitrarily. 
Let $p_n$ and $\tilde{p}_n$ denote the joint densities of $n$ consecutive observations from two zero-mean stationary processes, possibly non-Gaussian, with spectral densities $S$ and $\tilde S$, respectively. 
Without loss of generality, we normalize the time scale so that observations are taken at unit time intervals.

For $\alpha\in(0,1)$, a spectral $\alpha$-R\'{e}nyi divergence is defined as
\begin{align}
    &D_{\alpha}[\,S \,:\, \tilde{S}\,]
    := \frac{1}{2\pi(1-\alpha)}\int_{-\pi}^{\pi}
    [\log \{\alpha \tilde{S}(\omega)+(1-\alpha)S(\omega)\}
    \nonumber\\&\qquad\qquad\qquad\qquad\qquad\qquad\qquad-
    \alpha \log \tilde{S}(\omega)- (1-\alpha) \log S(\omega)] \intd \omega.
\end{align}
For $\alpha=1$, we define $D_{1}[\,S\,:\,\tilde{S}\,]$ by continuity as the limit $\alpha\to1^{-}$. 
In this limit, the spectral $\alpha$-R\'{e}nyi divergence converges to the Itakura--Saito divergence \cite{Pinsker1964,ItakuraSaito1968}:
\begin{align}
    \lim_{\alpha\to1^{-}}
    D_{\alpha}[\,S \,:\, \tilde{S}\,]
    &=D_{\mathrm{IS}}[\,S \,:\, \tilde{S}\,]
    \nonumber\\
    &:= \frac{1}{2\pi} \int_{-\pi}^{\pi}\left\{ 
    \left(\frac{\tilde{S}(\omega)}{S(\omega)}\right)^{-1}-1  
    +\log\left(\frac{\tilde{S}(\omega)}{S(\omega)} \right)
    \right\}\intd \omega.
\end{align}

The name of this class stems from the fact that the class is induced by the limit of the probabilistic R\'{e}nyi divergence between stationary Gaussian processes, as pointed out by \cite{Vajda1989}.
More precisely, if $p_n$ and $\tilde{p}_n$ are the joint densities of $n$ consecutive observations of two zero-mean stationary Gaussian processes with spectral densities $S$ and $\tilde{S}$, respectively, then for $\alpha\in(0,1]$, we have

    \[
    \lim_{n\to\infty}
    \frac{2}{n}
    D_{\alpha}[\,p_{n}\, :\, \tilde{p}_{n}\,]=
    D_{\alpha}[\,S\,:\,\tilde{S}\,],
    \]
    where $D_{\alpha}[\,p_{n}\, :\, \tilde{p}_{n}\,]$ is the (probabilistic) $\alpha$-R\'{e}nyi divergence:
    \[
    D_{\alpha}[\,p_{n}\, : \,\tilde{p}_{n}\,] := \frac{1}{\alpha-1}
    \log \int_{\mathbb{R}^{n}} \left(
    \frac{p_{n}(\mathbf{x}_{n})}{\tilde{p}_{n}(\mathbf{x}_{n})}
    \right)^{\alpha-1} p_{n}(\mathbf{x}_{n})\intd \mathbf{x}_{n}.
    \]
For $\alpha=1$, the probabilistic R\'{e}nyi divergence is defined by continuity  as the limit $\alpha\to1^{-}$ and coincides with the Kullback--Leibler divergence:
\begin{align*}
    D_{1}[\,p_n\,:\,\tilde{p}_n\,]
    = \int_{\mathbb{R}^n} p_n(\mathbf{x}_n) \log\frac{p_n(\mathbf{x}_n)}{\tilde{p}_n(\mathbf{x}_n)}\intd\mathbf{x}_n .
\end{align*}

The spectral $\alpha$-R\'{e}nyi divergence is actually a (statistical) divergence for any $\alpha\in(0,1)$:
that is, $D_{\alpha}[\,S\,:\,\tilde{S}\,]\ge 0$ and $D_{\alpha}[\,S\,:\,\tilde{S}\,]=0$ 
if and only if $S=\tilde{S}$.
This follows since the one-step ahead prediction error variance (the innovation variance) derived from the Kolmogorov–Szeg\"{o} formula (Theorem 5.8.1. of \cite{BrockwellandDavis})
\[
\sigma^{2}_{S}:=\exp\left(
    \frac{1}{2\pi}\int_{-\pi}^{\pi}\log{S(\omega)}\intd\omega
    \right)
\]
is log-concave with respect to $S$.

For $\alpha\in(0,1)$,
we shall prepare a discrete version of the spectral R\'{e}nyi divergence as follows:
\begin{align*}
&D_{\alpha}^{(n)}[\,S\,:\,\tilde{S}\,]
=
\frac{1}{(1-\alpha)n}
\sum_{\omega \in \Omega_{n}}
\Big{[}\log \left\{\alpha \tilde{S}\left(\omega\right)+ (1-\alpha)S\left(
\omega
\right)\right\}
\nonumber\\
&\qquad\qquad\qquad\qquad\qquad\qquad\qquad\qquad\qquad
 - \alpha \log \tilde{S}\left(
\omega
\right) - (1-\alpha) \log S\left(
\omega
\right)\Big{]},
\end{align*}
where $\Omega_{n}:=\{2\pi(t/n): t=-\lceil n/2\rceil+1,\ldots,-1,0,1,\ldots,\lfloor n/2 \rfloor \}$. 
By the convergence of the Riemann sums,
$D^{(n)}_{\alpha}$ converges to $D_{\alpha}$ as $n$ goes to infinity.
The discrete version $D^{(n)}_{\mathrm{IS}}$ of the Itakura--Saito divergence is defined in the same way.

\subsection{Properties of spectral R\'{e}nyi divergences}

Here we present two new formulae that have not been investigated in the literature:
(i) the asymptotic equivalence between the time-domain $\gamma$-divergence \cite{FujisawaEguchi2008} and the spectral R\'{e}nyi divergence; and
(ii) the variational representation of spectral R\'{e}nyi divergences.

We first state the asymptotic equivalence between the time-domain $\gamma$-divergence and the spectral R\'{e}nyi divergence.
Let $G_{\gamma}[\,p\,:\,q\,]$ be 
the $\gamma$-divergence between two probability densities $p$ and $q$:
\begin{align*}
&G_{\gamma}[\,p\,:\,q\,]
=\frac{1}{\gamma(1+\gamma)}
\log \int (p(x))^{1+\gamma}\intd x
-\frac{1}{\gamma}\log \int p(x)(q(x))^{\gamma}\intd x\\
&\qquad\qquad\qquad\qquad\qquad\qquad\qquad\qquad\qquad
+\frac{1}{1+\gamma}\log \int (q(x))^{1+\gamma}\intd x.
\end{align*}
\begin{theorem}[The time-domain $\gamma$-divergence leads to the spectral R\'{e}nyi divergence]
\label{thm: Gamma meets Renyi}
Let $\alpha\in(0,1)$ and $S,\tilde{S}\in\mathcal{S}$.
Let $p_n$ and $\tilde{p}_n$ denote the joint densities of $n$ consecutive observations from zero-mean stationary Gaussian processes with spectral densities $S$ and $\tilde{S}$, respectively. Then,
we have
\begin{align*}
\lim_{n\to\infty}\frac{2}{n}G_{\gamma}[\,p_{n}\,:\,\tilde{p}_{n}\,]
=D_{\alpha}[\, S \,:\, \tilde{S}\, ],
\end{align*}
where $\gamma=\alpha^{-1}-1>0$.
\end{theorem}

The proof is given in Appendix \ref{appendix: proof of Gamma}. It employs the direct evaluation of the $\gamma$-divergence and Szeg\"{o}'s limit theorem~\cite{Szego_1920,Kolmogorov_1941}.

We next state the variational representation of the spectral R\'{e}nyi divergence.
\begin{theorem}[Variational representation of the spectral R\'{e}nyi divergence]
\label{thm: variational rep e}
For $\alpha\in(0,1)$,
and
$S,\tilde{S}\in\mathcal{S}$,
we have
\begin{align}
D_{\alpha}
\,[\, S\,:\, \tilde{S} \,]
&=
\frac{1}{1-\alpha}
\min_{S'\in \mathcal{S}}
\left\{
\alpha \, D_{\mathrm{IS}}\,[\,\tilde{S}\,:\,S'\,]
+
(1-\alpha)\, D_{\mathrm{IS}}\,[\,S\,:\,S'\,]
\right\}
\label{eq: convex representation of Renyi e},
\end{align} 
where the minimum is achieved by $\underline{S}:=\alpha \tilde{S}+(1-\alpha)S$. The same characterization holds for $D^{(n)}_{\alpha}$.
\end{theorem}

\begin{proof}
Observe that for any $S'\in\mathcal{S}$, we have %\sout{$\alpha D_{\mathrm{IS}[\,\tilde{S}\,:\,S'\,]+ (1-\alpha)D_{\mathrm{IS}}[\,S\,:\,S'\,]}$}
\begin{align*}
& 
\alpha D_{\mathrm{IS}}[\,\tilde{S}\,:\,S'\,]
+ (1-\alpha)D_{\mathrm{IS}}[\,S\,:\,S'\,]
\\
&=\frac{1}{2\pi}\int_{-\pi}^{\pi}
\Bigg{\{}
\frac{\alpha \tilde{S}(\omega)+(1-\alpha)S(\omega)}{S'(\omega)}
-1+\log S'(\omega)\\
&\qquad\qquad\qquad\qquad\qquad\qquad\qquad
- \alpha \log \tilde{S}(\omega)
- (1-\alpha) \log S (\omega)
\Bigg{\}} \intd \omega \\
&=\frac{1}{2\pi}\int_{-\pi}^{\pi}
\Bigg{\{}
\frac{\alpha \tilde{S}(\omega)+(1-\alpha)S(\omega)}{S'(\omega)}
-1+\log \frac{S'(\omega)}{\alpha \tilde{S}(\omega)+(1-\alpha)S(\omega)}
\Bigg{\}} \intd \omega\\
& \,\, +\frac{1}{2\pi}
\int_{-\pi}^{\pi}[\log \{\alpha \tilde{S}(\omega)+(1-\alpha)S(\omega)\}-\alpha \log \tilde{S}(\omega)-(1-\alpha)\log S(\omega)]\intd \omega\\
&= D_{\mathrm{IS}}[\,\alpha \tilde{S}+(1-\alpha)S\,:\,S'\,] \\
& \,\, +\frac{1}{2\pi}
\int_{-\pi}^{\pi}[\log \{\alpha \tilde{S}(\omega)+(1-\alpha)S(\omega)\}-\alpha \log \tilde{S}(\omega)-(1-\alpha)\log S(\omega)]\intd \omega.
\end{align*}
So, the minimizer of the left hand side of the above equation is $\alpha \tilde{S}+(1-\alpha)S$, which proves the assertion.
\qed

The dual representation also holds:

\begin{proposition}
For $\alpha\in(0,1)$,
and
$S,\tilde{S}\in\mathcal{S}$,
we have
\begin{align}
D_{\alpha}
\,[\, S\,:\, \tilde{S} \,]&=
\frac{1}{1-\alpha}\min_{S'\in \mathcal{S}}
\left\{
\alpha \, D_{\mathrm{IS}}\,[\,S'\,:\,S\,]
+
(1-\alpha)\, D_{\mathrm{IS}}\,[\,S'\,:\,\tilde{S}\,]
\right\},
\label{eq: convex representation of Renyi m}
\end{align}
where the minimum is attained by $\overline{S}:=(\alpha S^{-1}+ (1-\alpha)\tilde{S}^{-1})^{-1}$.
The same characterization holds for $D^{(n)}_{\alpha}$.
\end{proposition}

\begin{proof}
Observe that we have, for any $S'\in\mathcal{S}$,
\begin{align*}
& \alpha D_{\mathrm{IS}}[\,S'\,:\,S\,]
+(1-\alpha) D_{\mathrm{IS}}[\,S'\,:\,\tilde{S}\,]\\
& = \frac{1}{2\pi}\int_{-\pi}^{\pi} \left[\left(\frac{\overline{S}(\omega)}{S'(\omega)}\right)^{-1}-1+\alpha \log S(\omega) + (1-\alpha) \log \tilde{S}(\omega) -\log S'(\omega)\right] \intd \omega\\
& = \frac{1}{2\pi}\int_{-\pi}^{\pi} \left[\left(\frac{\overline{S}(\omega)}{S'(\omega)}\right)^{-1}-1+ \log \frac{\overline{S}(\omega)}{S'(\omega)}\right] \intd \omega\\
& \qquad + \frac{1}{2\pi}\int_{-\pi}^{\pi}\left[
\alpha \log S(\omega)+(1-\alpha)\log \tilde{S}(\omega)-\log \overline{S}(\omega)
\right]\intd \omega\\
& = D_{\mathrm{IS}}[\,S'\,:\,\overline{S}\,]
+(1-\alpha)D_{\alpha}\,[\, S\,:\, \tilde{S} \,].
\end{align*}
So, the minimizer of the left hand side of the above equation is $\overline{S}$, which proves (\ref{eq: convex representation of Renyi m}).
\end{proof}

Note that this variational representation is a spectral version of the characterization of R\'{e}nyi divergence in relation to the composite hypothesis testing \cite{Shayevitz2010,vanErvenHarremos}.

\subsection{Robustness of the spectral R\'{e}nyi divergence}

A tantalizing feature of the spectral $\alpha$-R\'{e}nyi divergence 
is that it is robust with respect to outliers in the frequency domain, which has been discussed in \cite{ZhangTaniguchi1995,Hirukawa2005}.
Theorem \ref{thm: Gamma meets Renyi} explains this robustness as the time-domain $\gamma$-divergence induces the robust parameter estimation scheme \cite{FujisawaEguchi2008}.
Theorem \ref{thm: variational rep e} also provides another explanation about the robustness of the spectral R\'{e}nyi divergence as discussed in this subsection.
Let $\mathcal{S}_{\Theta}:=\{S_{\theta}\in\mathcal{S}: \theta\in\Theta\}$ be a parametric spectral model with a parameter space $\Theta\subset\mathbb{R}^{d}$, 
%\revisebegin
where $\Theta$ is compact with nonempty interior. 
%\reviseend

Consider the minimum spectral $\alpha$-R\'{e}nyi divergence estimator $\hat{\theta}_{\alpha}$ given as 
\begin{align*}
D^{(n)}_{\alpha}[\,\tilde{I}_{n}\,:\,S_{\hat{\theta}_{\alpha}}\,]
=\min_{\theta\in\Theta}D^{(n)}_{\alpha}[\,\tilde{I}_{n}\,:\,S_{\theta}\,],
\end{align*}
where $\tilde{I}_{n}$ is a nonparametric pilot estimate of the spectral density such as the periodogram
\begin{align*}
I_{n}(\omega):=\frac{1}{2\pi n}\left| \sum_{t=1}^{n}x_{t} \mathrm{e}^{-\sqrt{-1} t \omega} \right|^{2}
\end{align*}
and its smoothed versions.

We begin with introducing outliers and contamination in the frequency domain; see Section 8 of \cite{Maronna2019robust} for details.
Expressing contamination caused by outliers in the frequency domain,
we consider a contaminated pilot estimate $\tilde{I}_{n}^{z,\omega^{*}}$:
for $\omega^{*}\in\Omega_{n}$ and $z>0$,
\begin{align*}
\tilde{I}_{n}^{z,\omega^{*}}(\omega)=
\begin{cases}
\tilde{I}_{n}(\omega) + z,  & \omega=\pm \omega^{*},\\
\tilde{I}_{n}(\omega),      & \omega\ne \pm \omega^{*}.
\end{cases}
\end{align*}
This contamination is motivated by the case when
periodic components are not suitably subtracted and are contaminated in observed time series \cite{HeydeandDai1996,Iacone2010,McCloskeyandPerron2013}, say, $x_{t}^{o}=x_{t}+ \sqrt{8\pi z / n } \cos(At)$, $t=1,\ldots,n$ with original zero-mean time series $\{x_{t}\}_{t=1,2,\cdots,n}$. In this case, the periodogram of $x_{t}^{o}$ has approximately the form of 
\begin{align*}
I_{n}^{z,A}(\omega)=I_{n}(\omega)+z
(\delta_{A}(\omega)+\delta_{-A}(\omega)) + O_{P}(\sqrt{z}),
\end{align*}
where $\delta_{A}(\omega)=1$ if $\omega=A$ and 0 otherwise.

Observe first that the variation of the Itakura--Saito divergence with respect to an outlier in the frequency domain is written as
\begin{align}
D_{\mathrm{IS}}^{(n)}[\tilde{I}^{z,\omega^{*}}_{n}\,:\,S_{\theta}]
-
D_{\mathrm{IS}}^{(n)}[\tilde{I}_{n}\,:\,S_{\theta}]
=\frac{z}{n}\frac{1}{S_{\theta}(\omega^{*})}
-\frac{\log z}{n} + o_{P}(1)
\label{eq: change Itakura Saito}
\end{align}
with $o_{P}(1)$ being uniform with respect to $z$.
For larger value of $z$, this change is influenced by $\theta$, impacting on the behavior of the minimum Itakura--Saito divergence estimator.
Now,
putting the convex combination of $S_{\theta}$ and $\tilde{I}^{z,\omega^{*}}_{n}$ into
the second slot of the Itakura--Saito divergence mitigates the effect of $z$ as
\begin{align*}
&D_{\mathrm{IS}}^{(n)}[\tilde{I}^{z,\omega^{*}}_{n}\,:\,\alpha S_{\theta}+(1-\alpha)\tilde{I}^{z,\omega^{*}}_{n}]-
D_{\mathrm{IS}}^{(n)}[\tilde{I}_{n}\,:\, \alpha S _{\theta}+(1-\alpha)\tilde{I}_{n}]\\
&
= 
\frac{1}{n}\frac{\alpha S_{\theta}(\omega^{*})}{\alpha S_{\theta}(\omega^{*})+(1-\alpha)\tilde{I}_{n}(\omega^{*})}
\frac{z}{(1-\alpha)z+\alpha S_{\theta}(\omega^{*})+(1-\alpha)\tilde{I}_{n}(\omega^{*})}\\
&\quad-\frac{1}{n}\log (1+\tilde{I}_{n}^{-1}(\omega^{*})z)
 +\frac{1}{n}\log (1+(1-\alpha)z \{\alpha S_{\theta}(\omega^{*})+(1-\alpha)\tilde{I}_{n}(\omega^{*})\}^{-1})\\
&= o_{P}(1)
\end{align*}
with $o_{P}(1)$ being uniform with respect to $z$.
This, together with the variational representation 
\eqref{eq: convex representation of Renyi e}, implies that
\begin{align}
D_{\alpha}^{(n)}[\tilde{I}^{z,\omega^{*}}_{n}\,:\,S_{\theta}]
-D_{\alpha}^{(n)}[\tilde{I}_{n}\,:\, S_{\theta}]
=\frac{\alpha \log z}{n} + o_{P}(1)
\label{eq: change Renyi}
\end{align}
with $o_{P}(1)$ being uniform with respect to $z$.
So, for large $z$, the change does not vary with respect to $\theta$. This imbues the minimum spectral R\'{e}nyi divergence estimator with robustness against outliers in the frequency domain.

\section{Optimization paths}
\label{sec:outlier}

This section presents a further robustness property of the minimum spectral R\'{e}nyi divergence estimator. 
We show that
(1) the optimization for the minimum spectral R\'{e}nyi divergence estimator is robust against outliers in the frequency domain, whereas (2)
the optimization for
the minimum Itakura--Saito divergence estimator is sensitive to them.

\subsection{Stable optimization path of the Spectral R\'{e}nyi divergence minimization}

We begin with showing that the optimization for the minimum spectral R\'{e}nyi divergence estimator is robust against outliers in the frequency domain.
Consider the gradient descent update 
with a sequence of learning rates $\{\gamma_{i}>0:i=1,2,\ldots\}$ of $\hat{\theta}_{\alpha}$: for the $(m+1)$-th step,
\begin{align*}
\hat{\theta}^{(m+1)}_{\alpha}[\tilde{I}_{n}]=
\hat{\theta}^{(m)}_{\alpha}[\tilde{I}_{n}]
+ \gamma_{m}\,
\mathcal{G}_{\alpha}(\hat{\theta}^{(m)}_{\alpha}[\tilde{I}_{n}]\,;\,\tilde{I}_{n})
\end{align*}
with
\begin{align}
\mathcal{G}_{\alpha}(\theta\,;\,\tilde{I}_{n})
&=
-\nabla_{\theta}D^{(n)}_{\alpha}[\tilde{I}_{n}\,:\,S_{\theta}]
\nonumber\\
&=
\frac{\alpha}{(1-\alpha) n }\sum_{\omega\in\Omega_{n}}
\left[1-\frac{S_{\theta}(\omega)}{\alpha S_{\theta}(\omega)+(1-\alpha)\tilde{I}_{n}(\omega)}\right]
 \nabla_{\theta}\log S_{\theta}(\omega).
 \label{expr:grad-SpecRenyi}
\end{align}

For theoretical results, we make the following assumptions. 
%\revisebegin 
Let $\|\cdot\|$ denote the Euclidean norm. 
%\reviseend
\begin{assumption}\rm\label{Assumption: norm constraint of spectral}
    %\revisebegin
    Let $S_{\theta}\in\mathcal{S}_{\Theta}$. For each $\omega \in (0,\pi)$, the function $\theta\mapsto S_\theta(\omega)$ is twice continuously differentiable on $\mathrm{int}(\Theta)$, where $\mathrm{int}(\Theta)$ denotes the interior of $\Theta$. 
    Moreover, the following conditions hold:
    %\reviseend
    \begin{align*}
        &\sup_{\theta\in \mathrm{int}(\Theta),\,\omega\in[-\pi,\pi]}\max\{\|\nabla_{\theta} S_{\theta}(\omega)\|,\,
        \|\nabla_{\theta}\log S_{\theta}(\omega)\|
        ,\,
        \| (1/S_{\theta}(\omega))\nabla_{\theta} \log S_{\theta}(\omega)\|\}\\
        &\qquad\qquad\qquad=: U_{1} <\infty; %\\
    \end{align*}
    and
    \begin{align*}
        %&\qquad\qquad 
        \limsup_{n\to\infty}
        \left(\sup_{\theta\in\mathrm{int}(\Theta)}\frac{1}{n}\sum_{\omega\in\Omega_{n}}\|\nabla^{2}_{\theta}\log S_{\theta}(\omega)\|_{\mathrm{op}} \right) =: U_{2} <\infty,
    \end{align*}
    where $\|\cdot\|_{\mathrm{op}}$ denotes the operator norm.
    %\reviseend
\end{assumption}

\begin{assumption}\rm
\label{Assumption: bounded path}
 For any $z\ge 0$, $\omega^{*}\in\Omega_{n}$, and $m\in\mathbb{N}$,
the $m$-th gradient descent update $\hat{\theta}^{(m)}_{\alpha}[\tilde{I}^{z,\omega^{*}}]$
lies in %\sout{$\Theta$} 
%\revisebegin
$\mathrm{int}(\Theta)$. 
%\reviseend
\end{assumption}

Under Assumptions \ref{Assumption: norm constraint of spectral} and \ref{Assumption: bounded path},
the subsequent theorem indicates that 
for any value of $\alpha\in(0,1)$,
the path of the gradient descent 
with a prespecified step size
for $\hat{\theta}_{\alpha}$ 
is robust against the presence of outliers in the frequency domain.

\begin{theorem}[Stability of optimization path with a fixed learning rate sequence]
\label{thm: success Renyi}
Fix arbitrary step size sequence $\{\gamma_{k}>0\,:\,k=1,2,\ldots\}$.
Under Assumptions 
$\ref{Assumption: norm constraint of spectral}$ and $\ref{Assumption: bounded path}$,
for any $m\in\mathbb{N}$, we have
\begin{align}
\lim_{n\to\infty}
\sup_{z\ge 0,\omega^{*}\in\Omega_{n}}
\left\|
\hat{\theta}_{\alpha}^{(m)}[\tilde{I}^{z,\omega^{*}}_{n}]
-
\hat{\theta}_{\alpha}^{(m)}[\tilde{I}_{n}]
\right\|
=0 \text{ almost surely}.
\end{align}
\end{theorem}

The proof is given in Appendix \ref{appendix: proof of thm}.

\begin{remark}\rm(Constant step size)~
Assumption $\ref{Assumption: norm constraint of spectral}$ ensures 
that we can obtain an upper bound of the smoothness constant of $D^{(n)}_{\alpha}[\tilde{I}_{n}^{z,\omega^{*}}:S_{\theta}]$ with respect to $\theta$
uniformly in $z\ge 0$: for any $z\ge 0$, we have
\begin{align}
&\sup_{\theta\ne\theta'}
\frac{
\|\nabla_{\theta}D^{(n)}_{\alpha}[\tilde{I}^{z,\omega^{*}}_{n}\,:\,S_{\theta}]
-
\nabla_{\theta}D^{(n)}_{\alpha}
[\tilde{I}^{z,\omega^{*}}_{n}\,:\,S_{\theta'}]
\|}{\|\theta-\theta'\|}\nonumber\\
&\le \sup_{z\ge 0}\sup_{\theta\in 
%\mbox{\sout{$\Theta$}}
%\revisebegin
\mathrm{int}(\Theta). 
%\reviseend
}\|\nabla_{\theta}^{2}D_{\alpha}^{(n)}[\tilde{I}_{n}^{z,\omega^{*}}:S_{\theta}]\|_{\mathrm{op}}\nonumber\\
&\le L:= \frac{U_{1}^{2}+U_{2}}{1-\alpha}<\infty,
\label{eq: Lipschitz}
\end{align}
where we use the following expression of the Hessian of $D_{\alpha}^{(n)}[\tilde{I}_{n}:S_{\theta}]$:
\begin{align*}
&\left(-\nabla_{\theta}^{2}
D^{(n)}_{\alpha}[\tilde{I}_{n}:S_{\theta}]\right)\\
&=-\frac{1}{(1-\alpha)n}\sum_{\omega\in\Omega_{n}}\frac{\alpha S_{\theta}(\omega)}
{\alpha S_{\theta}(\omega) + (1-\alpha)\tilde{I}_{n}(\omega)}\\
&\qquad\qquad\qquad\qquad\frac{(1-\alpha) \tilde{I}_{n}(\omega)}
{\alpha S_{\theta}(\omega) + (1-\alpha)\tilde{I}_{n}(\omega)}
\left(\nabla_{\theta}\log S_{\theta}(\omega)\right)
\left(\nabla_{\theta}\log S_{\theta}(\omega)\right)^{\top}\\
&\quad 
+\frac{\alpha}{(1-\alpha)n}\sum_{\omega\in\Omega_{n}}
\left\{1-
\frac{S_{\theta}(\omega)}
{\alpha S_{\theta}(\omega)+(1-\alpha)\tilde{I}_{n}(\omega)}\right\}
\nabla^{2}_{\theta}\log S_{\theta}(\omega).
\end{align*}
By the standard theory of the gradient descent algorithm (c.f.,~Theorem $3.2$ of \cite{Nocedal_Wright_book}),
for any $z\ge 0$,
the gradient along with
the gradient descent update 
$\{\hat{\theta}^{(m)}_{\alpha}[\tilde{I}^{z,\omega^{*}}]\}_{m=0,1,2,\ldots}$
with a constant step size $\gamma_{m}= L^{-1}$
vanishes: 
\[
\lim_{m\to\infty}
\nabla_{\theta=\hat{\theta}^{(m)}_{\alpha}[\tilde{I}_{n}^{z,\omega^{*}}]}
D_{\alpha}^{(n)}[\tilde{I}_{n}^{z,\omega^{*}}
:
S_{\theta}]
=0.
\]
If the further assumption called the Kurdyka–{\L}ojasiewicz inequality \cite{Lojasiewicz_1963,Kurdyka_1998} holds 
at all stationary points, that is, points satisfying
$\nabla_{\theta}
D_{\alpha}^{(n)}[\tilde{I}_{n}^{z,\omega^{*}}
:
S_{\theta}]=0$
for the function $\theta\mapsto D_{\alpha}^{(n)}[\tilde{I}_{n}^{z,\omega^{*}}:S_{\theta}]$,
the bounded gradient descent update converges to a stationary point;
c.f., Theorem $3.2$ of \cite{AttouchBolteSvaiter2013}.
Theorem $\ref{thm: success Renyi}$
implies that the finite update sequences
for different $z\ge 0$
become the same 
and all stationary points for different $z\ge 0$ become the same
as $n\to\infty$, which provides a new perspective on robust spectral analysis and optimization theory.
\end{remark}

\begin{remark}\rm(AR(1) model)~
    Here we check Assumption $\ref{Assumption: norm constraint of spectral}$ for a first-order autoregressive model. The model spectral density is
    \[
    S_{\sigma,\rho}(\omega) = \frac{1}{2\pi}\frac{\sigma^{2}}{1-2\rho \cos(\omega)+\rho^{2}},\quad \sigma>0\quad \text{and}\quad \rho\in(-1,1).
    \]
    We reparametrize this by using 
    \[
    \theta_{1} = \log \sigma \quad \text{and} \quad \theta_{2} = 
    \,\mathrm{atanh}(\rho)
    \]
    %\revisebegin
    with 
    \[
    \mathrm{atanh}(\rho)=\frac{1}{2}\log\frac{1+\rho}{1-\rho}
    \]
    %\reviseend
    so as to remove constraints in the optimization.
    Then, the gradient of $\log S_{\theta}(\omega)$ with respect to $\theta$ is given by 
    \[
    \nabla_{\theta} \log S_{\theta} (\omega) = \begin{pmatrix}
        2 \\
        2(1-\rho^{2})\{\cos(\omega)-\rho\}/\{1-2\rho \cos(\omega)+\rho^{2}\}
    \end{pmatrix}.
    \]
    The Hessian of $\log S_{\theta}(\omega)$ is a bit complex but is bounded
    as
    \[
    \| \nabla^{2}_{\theta} \log S_{\theta}(\omega) \|_{\mathrm{op}}
    \le 
    C / (1-|\rho|)^{4}
    \]
    for some absolute positive constant $C$.
    %\sout{So, Assumption $\ref{Assumption: norm constraint of spectral}$ holds for any bounded parameter space $\Theta$ of $\theta$.}
    %\revisebegin 
    So, Assumption $\ref{Assumption: norm constraint of spectral}$ holds on the interior of any compact subset $\Theta \subset (0,\infty)\times(-1,1)$.
    %\reviseend
\end{remark}

\begin{figure}[htbp]
    \centering
\includegraphics[scale=0.43]{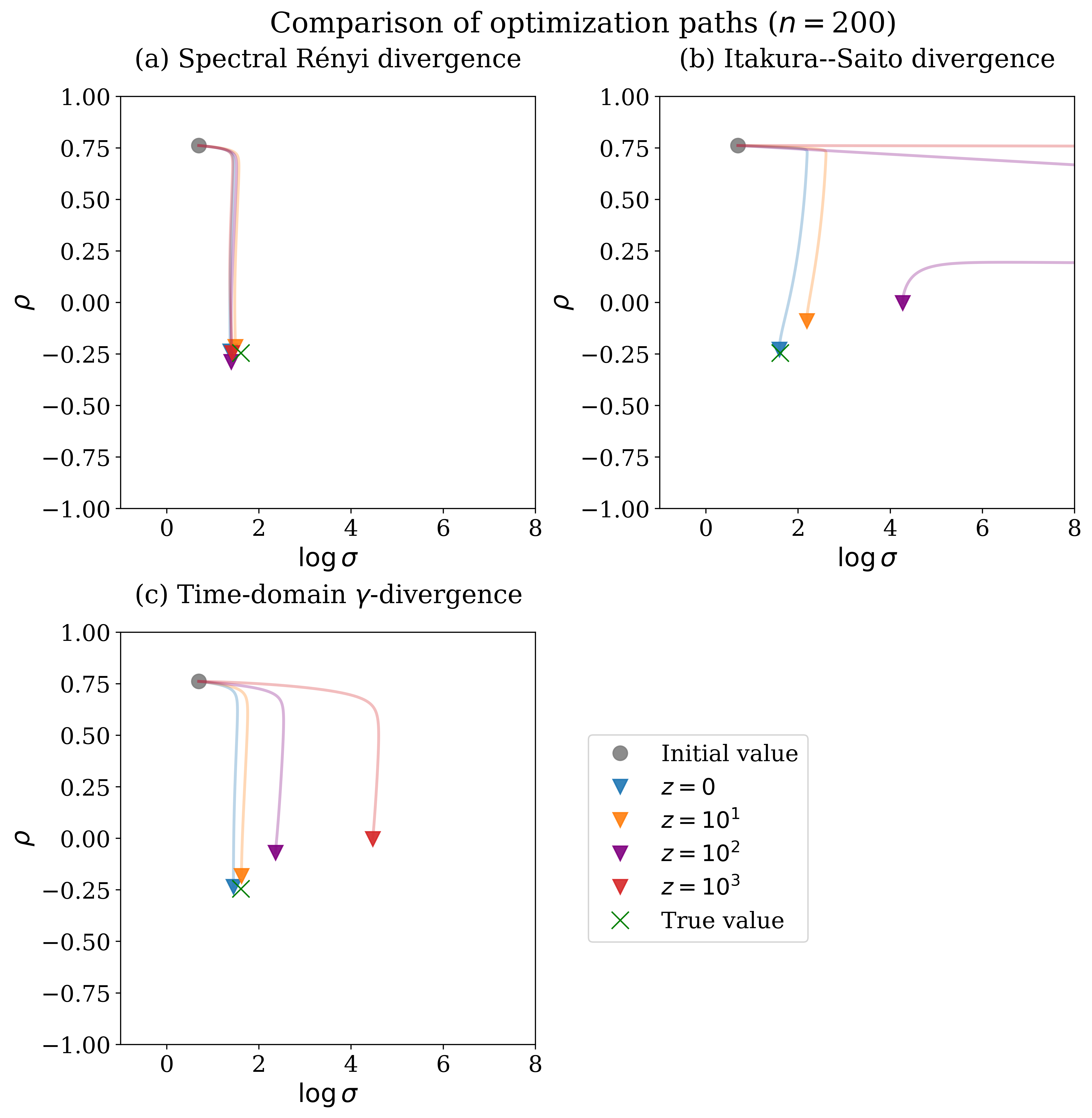}
    \caption{Comparison of optimization paths when the length $n$ of a time series is $200$.
    The optimization paths of the gradient descent for the 
    %\revisebegin
        AR(1) 
    %\reviseend
    %\sout{$AR(1)$} 
    model, with respect to $(\log \sigma, (1-\exp(-2\rho))/(1+\exp(-2\rho)))$, are presented under the spectral R\'{e}nyi divergence minimization, the Itakura--Saito divergence minimization, and the time-domain $\gamma$-divergence minimization.
    Panel (a) shows the optimization paths (curves with light colors) of the spectral R\'{e}nyi divergence minimization that start from the initial value (the gray circle) and terminate at points after $10000$ iterations (the inverted triangles), where the colors indicate different values of the contamination $z$.
    Panels (b) and (c) show the corresponding results for the Itakura--Saito divergence  and time-domain $\gamma$-divergence minimization, respectively. 
    }
    \label{Fig_Renyi_demonstrate}
\end{figure}

Theorem \ref{thm: success Renyi} establishes a connection between robustness in statistics and optimization theory, offering a novel perspective in robust statistics. While detailed numerical experiments are presented in Section \ref{sec:numerical}, we briefly illustrate the implications through a simple simulation study here.

We consider contaminated observations defined as
\[
X_{t}^{\circ}=X_{t}+\sqrt{z}\sin(t\pi/2), \quad t=1,\ldots,n=200
\]
with $X_{t}$ generated from AR(1) model having $(\sigma,\rho)$.
We perform 10000 iterations of gradient descent with a step size of 0.01
with respect to \[ ( \log \sigma, 
(1-\exp(-2\rho))/(1+\exp(-2\rho)) ). \]
Figure \ref{Fig_Renyi_demonstrate} (a) illustrates optimization paths and endpoints based on \(D^{(n)}_{\alpha=0.5}[I_{n} : S_{\theta}]\), starting from the initial value (gray circle). Different colors represent different contamination levels \(z \in \{0, 10^{1}, 10^{2}, 10^{3}\}\).
The figure highlights the robustness of the optimization paths under varying contamination levels.
In contrast, Figure \ref{Fig_Renyi_demonstrate} (b) shows the optimization paths and the stopping points based on the Itakura--Saito divergence 
$D^{(n)}_{\mathrm{IS}}[I_{n}\,:\,S_{\theta}]$ starting from the initial value (the gray circle).
This figure demonstrates the non-robustness of the optimization paths of the Itakura--Saito divergence minimization under contamination. The instability observed here will be analyzed further in the subsequent subsection.
Further, Figure \ref{Fig_Renyi_demonstrate} (c) shows the optimization paths and the endpoints based on the time-domain $\gamma$-divergence. Although Theorem~\ref{thm: Gamma meets Renyi} establishes a link between the time-domain $\gamma$-divergence and the spectral $\alpha$-R\'{e}nyi divergence, their non-asymptotic behavior is different and
the implementation of the time-domain $\gamma$-divergence minimization
becomes increasingly difficult as $n$ grows, because each optimization step requires a costly log-determinant computation; see also the proof of Theorem~\ref{thm: Gamma meets Renyi} in the supplement. This comparison highlights a practical advantage of the approach using the spectral R\'{e}nyi divergence.

Practically, gradient descent updates commonly use line search methods to optimize the step size dynamically. 
To manage this situation, we present a guiding theorem for gradient descent updates with line search. In this paper, we employ the Armijo condition for this purpose.
Consider an objective function $\mathcal{L}(\theta)$ to be minimized and let $\theta^{(m)}$ be the $m$-th gradient descent update for each $m\in\mathbb{N}$.
The Armijo condition with $c\in(0,1)$ selects the step size $\gamma$ in the $m$-th update $\theta^{(m)}$ so as to satisfy
\[
\mathcal{L}\left(\theta^{(m-1)} - \gamma \nabla_{\theta=\theta^{(m-1)}} \mathcal{L}(\theta) \right) 
\le 
\mathcal{L}\left(\theta^{(m-1)}\right) 
-c \gamma \left\|\nabla_{\theta=\theta^{(m-1)}}\mathcal{L}(\theta) \right\|^{2},
\]
which ensures the monotone decrease of the objective function along with the gradient descent update.
We denote by $\mathcal{A}[c,\theta,\mathcal{L}]$ the Armijo condition with $c\in(0,1)$ for an objective function $\mathcal{L}(\theta)$ in the update from $\theta$.

\begin{proposition}[Stability of line search]
\label{prop: stability of line search}
Fix $c\in(0,1)$, $\kappa>0$, and $\underline{\gamma}>0$.
Fix also the point $\theta\in\Theta$ such that
\[
\left\|\nabla_{\theta}D_{\alpha}^{(n)}[\tilde{I}_{n}:
S_{\theta}] \right\|>\kappa.\]
Then, 
under Assumption $\ref{Assumption: norm constraint of spectral}$,
for sufficiently large $n$,
the step size $\gamma$ satisfying 
$\mathcal{A}[c\,,\,\theta\,,\,D_{\alpha}^{(n)}[\tilde{I}_{n} :S_{\theta}]]$
and $\gamma\in (\underline{\gamma}, 2(1-c)/L)$
also satisfies $\mathcal{A}[\tilde{c}\,,\,\theta\,,\,D_{\alpha}^{(n)}[\tilde{I}_{n}^{z,\omega^{*}} :S_{\theta}]]$
for any $\tilde{c}\in(0,c)$, $z>0$, and $\omega^{*}\in\Omega_{n}$,
where $L$ is the smoothness constant in $(\ref{eq: Lipschitz})$.
\end{proposition}

The proof is given in Appendix \ref{appendix: Proof of line search}.
It employs the descent lemma (c.f., \cite{Bertsekas_NP}) together with the comparison of gradients.
\begin{remark}\rm(Optimization path with the Armijo condition)~
Let $\{\gamma_{m}[\tilde{I}_{n}]:m=1,2,\ldots\}$
and $\{\gamma_{m}[\tilde{I}^{z,\omega^{*}}_{n}]:m=1,2,\ldots\}$ 
the step size sequences for $D_{\alpha}^{(n)}[\tilde{I}_{n}:S_{\theta}]$
and $D_{\alpha}^{(n)}[\tilde{I}_{n}^{z,\omega^{*}}:S_{\theta}]$, respectively.
From Proposition $\ref{prop: stability of line search}$, we can assume
\[
\sup_{z\ge 0,\omega^{*}\in\Omega_{n}}|\gamma_{m}[\tilde{I}_{n}]-\gamma_{m}[\tilde{I}_{n}^{z,\omega^{*}}]|\to 0 \quad \text{as }\quad n\to\infty.
\]
Under this additional assumption and 
if we have
\[\inf_{k=0,\ldots,m-1}\|\nabla_{\theta=\hat{\theta}_{\alpha}^{(k)}[\tilde{I}_{n}]}D_{\alpha}^{(n)}[\tilde{I}_{n}:
S_{\theta}] \|>0, \]
we can have
for any $m\in\mathbb{N}$,
\begin{align*}
\lim_{n\to\infty}
\sup_{z\ge 0,\omega^{*}\in\Omega_{n}}
\left\|
\hat{\theta}_{\alpha}^{(m)}[\tilde{I}^{z,\omega^{*}}_{n}]
-
\hat{\theta}_{\alpha}^{(m)}[\tilde{I}_{n}]
\right\|
=0 \text{ almost surely},
\end{align*}
even when we utilize the line search.
\end{remark}

\subsection{Unstable optimization path of the Itakura--Saito divergence minimization}

We next show that the optimization for
the minimum Itakura--Saito divergence estimator is sensitive to outliers in the frequency domain.
Consider the gradient descent update with a sequence of learning rates $\{\gamma_{i}>0:i=1,2,\ldots\}$
of $\hat{\theta}_{\mathrm{IS}}$: for the $(m+1)$-th step,
\begin{align*}
\hat{\theta}_{\mathrm{IS}}^{(m+1)}[\tilde{I}_{n}]
=\hat{\theta}_{\mathrm{IS}}^{(m)}[\tilde{I}_{n}] + \gamma_{m}
\, \mathcal{G}_{\mathrm{IS}}(\hat{\theta}_{\mathrm{IS}}^{(m)}[\tilde{I}_{n}]\,;\, \tilde{I}_{n})
\end{align*}
with
\begin{align*}
\mathcal{G}_{\mathrm{IS}}(\theta\,;\, \tilde{I}_{n})
=
\frac{1}{ n}
\sum_{\omega\in\Omega_{n}}
\left[
\frac{\tilde{I}_{n}(\omega)}{S_{\theta}(\omega)} - 1 \right]
\nabla_{\theta}
\log S_{\theta}(\omega).
\end{align*}
The subsequent proposition suggests that
both the initial update and the convergent point for $\hat{\theta}_{\mathrm{IS}}$ are highly sensitive to outliers in the frequency domain. %\sout{{Let $\|\cdot\|$ denote the Euclidean norm.}}
\begin{proposition}
\label{prop: failure Itakura Saito}
Fix $n\in\mathbb{N}$.
Assume that for $\omega^{*}\in\Omega_{n}$,
the inequality
\[
\inf_{\theta\in %\mbox{\sout{$\Theta$}}
%\revisebegin
\mathrm{int}(\Theta)
%\reviseend
} \|
(1/S_{\theta}(\omega^{*}))
\nabla_{\theta} \log S_{\theta}
(\omega^{*})\|>0
\] holds and $\theta^{(0)}$ lies in %\sout{$\Theta$}
%\revisebegin
$\mathrm{int}(\Theta)$. 
%\reviseend
Then we have
\begin{align}
\sup_{z>0}\left\|\hat{\theta}^{(1)}_{\mathrm{IS}}[\tilde{I}^{z,\omega^{*}}_{n}]
-
\hat{\theta}^{(1)}_{\mathrm{IS}}[\tilde{I}_{n}]
\right\|=\infty \,\,\text{almost surely}.
\end{align}
Further, 
assume that
for any $z>0$ and $\omega^{*}\in \Omega_{n}$,
the gradient descent sequence $\{\hat{\theta}^{(m)}_{\mathrm{IS}}[\tilde{I}^{z,\omega^{*}}]\}_{m=1,2\ldots}$
lies in %\sout{$\Theta$}
%\revisebegin
$\mathrm{int}(\Theta)$. 
%\reviseend
Then,
even when 
the gradient descent sequence with $\tilde{I}^{z,\omega^{*}}_{n}$ has a subsequence converging to a stationary point $\hat{\theta}^{(\infty)}_{\mathrm{IS}}[\tilde{I}^{z,\omega^{*}}_{n}]$ 
for $D_{\mathrm{IS}}^{(n)}[\,\tilde{I}^{z,\omega^{*}}_{n}\,:\,S_{\theta}\,]$, {\it i.e.} \[\nabla_{\theta=\hat{\theta}^{(\infty)}_{\mathrm{IS}}[\tilde{I}^{z,\omega^{*}}_{n}]}D_{\mathrm{IS}}^{(n)}[\tilde{I}_{n}^{z,\omega^{*}}\,:\,S_{\theta}]=0,\]
this stationary point $\hat{\theta}^{(\infty)}_{\mathrm{IS}}[\tilde{I}^{z,\omega^{*}}_{n}]$ is not a stationary point for $D_{\mathrm{IS}}^{(n)}[\,\tilde{I}_{n}\,:\,S_{\theta}\,]$:
\begin{align*}
\sup_{z>0}
\left\|\nabla_{\theta=\hat{\theta}^{(\infty)}[\tilde{I}_{n}^{z,\omega^{*}}]}D_{\mathrm{IS}}^{(n)}[\,\tilde{I}_{n}\,:\,S_{\theta}\,]\right\|=\infty
\,\,\text{almost surely}.
\end{align*}

\end{proposition}

The proof is given in Appendix \ref{appendix: proof of Prop}.

\section{Simulation studies}
\label{sec:numerical}

This section presents two numerical studies using the AR(2) model and the Brune spectral model with attenuation that is motivated by seismological studies \cite{Yoshimitsuetal2023}.

%\begin{comment}
\subsection{AR(2) model}

\begin{table}[ht]
\centering
\caption{
%\revisebegin 
The mean values of estimates with standard deviations for the AR(2) model without any trend. The values closest to the true values are underlined.
%\reviseend
}
\begin{tabular}{r|cccc}
  & $\hat{\sigma}$ & $\hat{\varphi_{1}}$ & $\hat{\varphi_{2}}$ \\ 
  \hline
  R\'{e}nyi with $\alpha=0.50$ & 
  0.83 ($\pm$ 0.06) & 0.57 ($\pm$ 0.64) & 0.14 ($\pm$ 0.65) \\ 
  R\'{e}nyi with $\alpha=0.75$ & 
  0.95 ($\pm$ 0.10) & \underline{0.84} ($\pm$ 0.50) & \underline{-0.17} ($\pm$ 0.44) \\ 
  R\'{e}nyi with $\alpha=0.90$ & 
  \underline{1.02} ($\pm$ 0.12) & 1.01 ($\pm$ 0.40) & -0.28 ($\pm$ 0.38) \\ 
  IS with the periodogram  & 
  1.13 ($\pm$ 0.18) & 1.09 ($\pm$ 0.48) & -0.32 ($\pm$ 0.42) \\ 
  IS with the smoothed periodogram & 
  1.14 ($\pm$ 0.19) & 1.11 ($\pm$ 0.51) & -0.34 ($\pm$ 0.44) \\ 
   \hline
\end{tabular}
\label{tab: AR without trend}
\end{table}

\begin{table}[ht]
\centering
\caption{
%\revisebegin
The mean values of estimates with standard deviations for the AR(2) model with a trigonometric trend. The values closest to the true values are underlined.
%\reviseend
} 
\begin{tabular}{r|cccc}
  & $\hat{\sigma}$ & $\hat{\varphi_{1}}$ & $\hat{\varphi_{2}}$ \\ 
  \hline
R\'{e}nyi with $\alpha=0.50$ & 0.85 ($\pm$ 0.05) & 0.53 ($\pm$ 0.68) & 0.20 ($\pm$ 0.74) \\ 
  R\'{e}nyi with $\alpha=0.75$ & \underline{0.99} ($\pm$ 0.12) &\underline{0.81} ($\pm$ 0.58) & \underline{-0.14} ($\pm$ 0.51) \\ 
  R\'{e}nyi with $\alpha=0.90$ & 1.07 ($\pm$ 0.14) & 1.02 ($\pm$ 0.46) & -0.26 ($\pm$ 0.42) \\ 
  IS with the periodogram  & 2.35 ($\pm$ 0.09) & 0.02 ($\pm$ 0.05) & 0.98 ($\pm$ 0.05) \\ 
  IS with the smoothed periodogram & 2.35 ($\pm$ 0.10) & 0.03 ($\pm$ 0.06) & 0.97 ($\pm$ 0.06) \\ 
   \hline
\end{tabular}
\label{tab: AR with trend}
\end{table}

First, we check the performance using the AR(2) model, where the spectral density is given by
\begin{align*}
    S^{\mathrm{AR2}}_{\theta}(\omega):=
    \sigma^{2} \Big{/} \left|1-\varphi_{1} e^{-\sqrt{-1}\omega} -\varphi_{2} e^{-2\sqrt{-1}\omega}\right|^{2},\;\; \theta=(\sigma,\varphi_{1},\varphi_{2}).
\end{align*}
We use the true value 
$\theta^\ast:=(\sigma^\ast,\varphi_{1}^{\ast},\varphi_{2}^{\ast})=(1,0.9,-0.2)$.
We use three different initial values:
$\theta^{(0),1}=(0.1,0,0)$,
$\theta^{(0),2}=(1,0,0)$,
and
$\theta^{(0),3}=(10,0,0)$.
%\revisebegin
The optimization is constrained to the stationary region of the AR(2) process; see Section~\ref{appendix: additional comparison} for details.
%\reviseend
We generate a time series $\mathbf{x}_{n=2^{10}}$ with its spectral density $S_{\theta^{*}}^{\mathrm{AR2}}$
1000 times,
and if we consider the contamination in the frequency domain,
we add
trigonometric trends to the time series as
\[
x_{t}^{o}=x_{t}+z_{1}\sqrt{2\pi/n}\sin(A_{1}t)+z_{2}\sqrt{2\pi/n}\sin(A_{2}t),\,\,
t=1,\ldots,n
\]
with $z_{1}=z_{2}=20$, $A_{1}=\pi/4$ and $A_{2}=\pi/8$.

We compare five estimation methods:
(1,2,3) the minimum spectral R\'{e}nyi divergence estimators with 
\[\alpha\in\{0.5,0.75,0.9\}\] using the 
%\revisebegin
%\sout{smoothed}
%\reviseend
periodogram ($I_{n}$);
(4) the minimum Itakura--Saito divergence estimator employing the  periodogram ($I_{n}$); and
(5) the minimum Itakura--Saito divergence estimator with the smoothed periodogram ($I_{n}^{S}$).

For each divergence minimization problem, the minimization is conducted via the gradient descent algorithm 
with a fixed learning rate $\gamma=0.005$. We stop the optimization when the optimization step reaches 10000 or when  
the Euclidean norm of the gradient of each divergence evaluated at each updated estimate is less than $10^{-3}$.
For smoothing the periodogram,
we apply the modified Daniell smoothers twice, where the length of the first time is $3$ and the length of the second time is $5$.

Tables \ref{tab: AR without trend} and \ref{tab: AR with trend} summarize the estimation results in the above setup. 
Tables~\ref{tab: AR without trend} and \ref{tab: AR with trend} suggest that the minimizers of the IS divergence with the periodogram/smoothed periodogram are very sensitive to the trigonometric trends of the observed time series~({\it i.e.} the outliers in the frequency domain). 
On the other hand, the minimizers of the spectral R\'enyi divergence give stably estimated values and their performance would be substantially better than those of the IS divergence with the periodogram/smoothed periodogram. 
Thus the estimated values using the spectral R\'enyi divergence are robust against the outliers in the frequency domain.

%\revisebegin
\begin{figure}[h]
    \centering
    \includegraphics[scale=0.37]{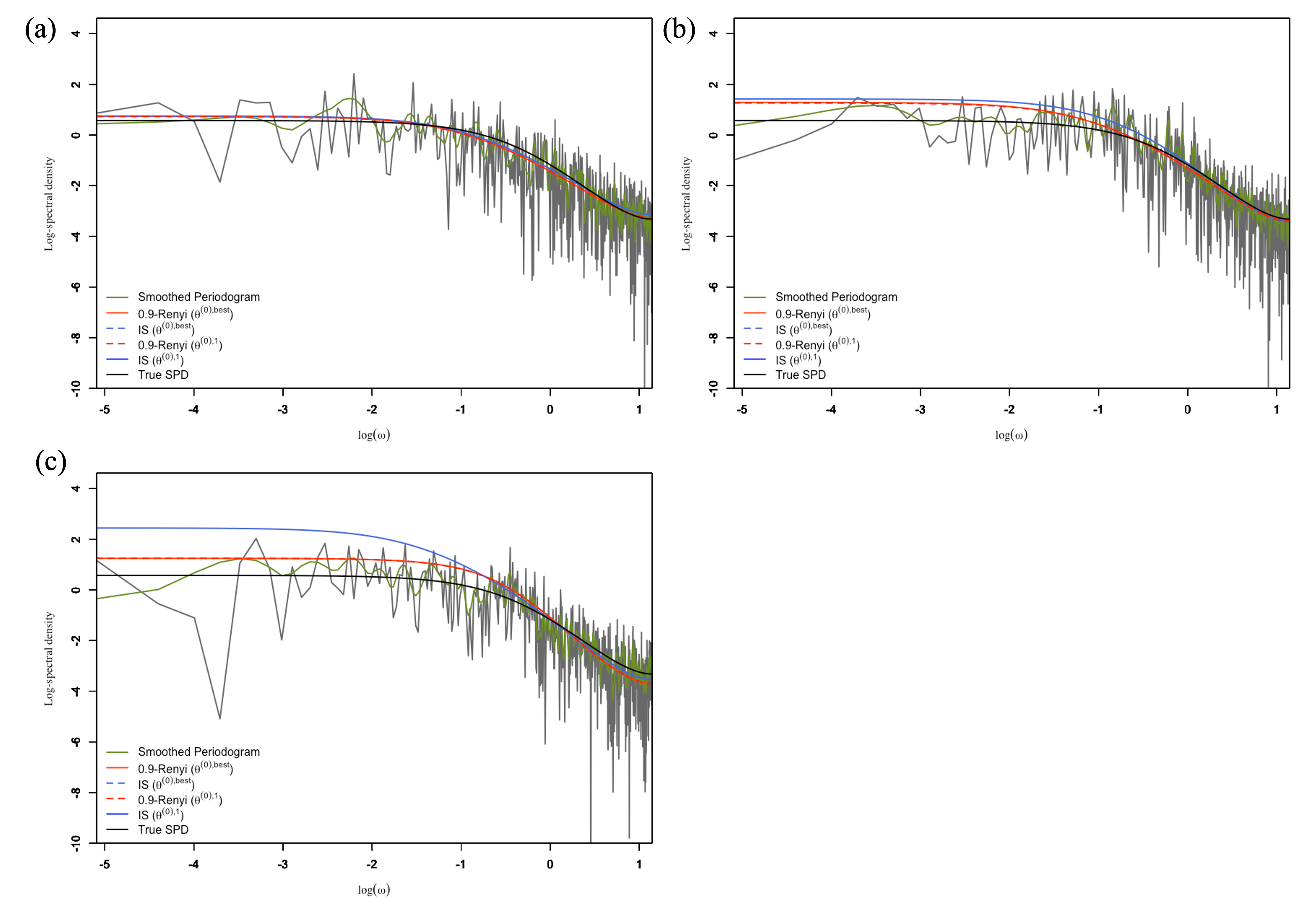}
    \caption{ %\revisebegin 
    Spectral 
    %\reviseend 
    densities with the estimates plugged in for the AR(2) model without any trend. 
    Each panel shows the results for different realizations of the periodogram. 
    The gray curve is the periodogram. 
    The green curve is the smoothed periodogram.
    The true spectral density is colored in black. The spectral density based on the spectral R\'{e}nyi divergence is colored in red. The spectral density based on the Itakura--Saito divergence is colored in blue.
    We set $\theta^{(0),1}$ to the initial value.
    The initial value $\theta^{(0),\mathrm{best}}$ represents the point that minimizes the objective functions over a grid of the parameter space. %, respectively.
    }
    \label{Fig_SPD1_AR2}
\end{figure}

%\revisebegin
\begin{figure}[h]
    \centering
    \includegraphics[scale=0.37]{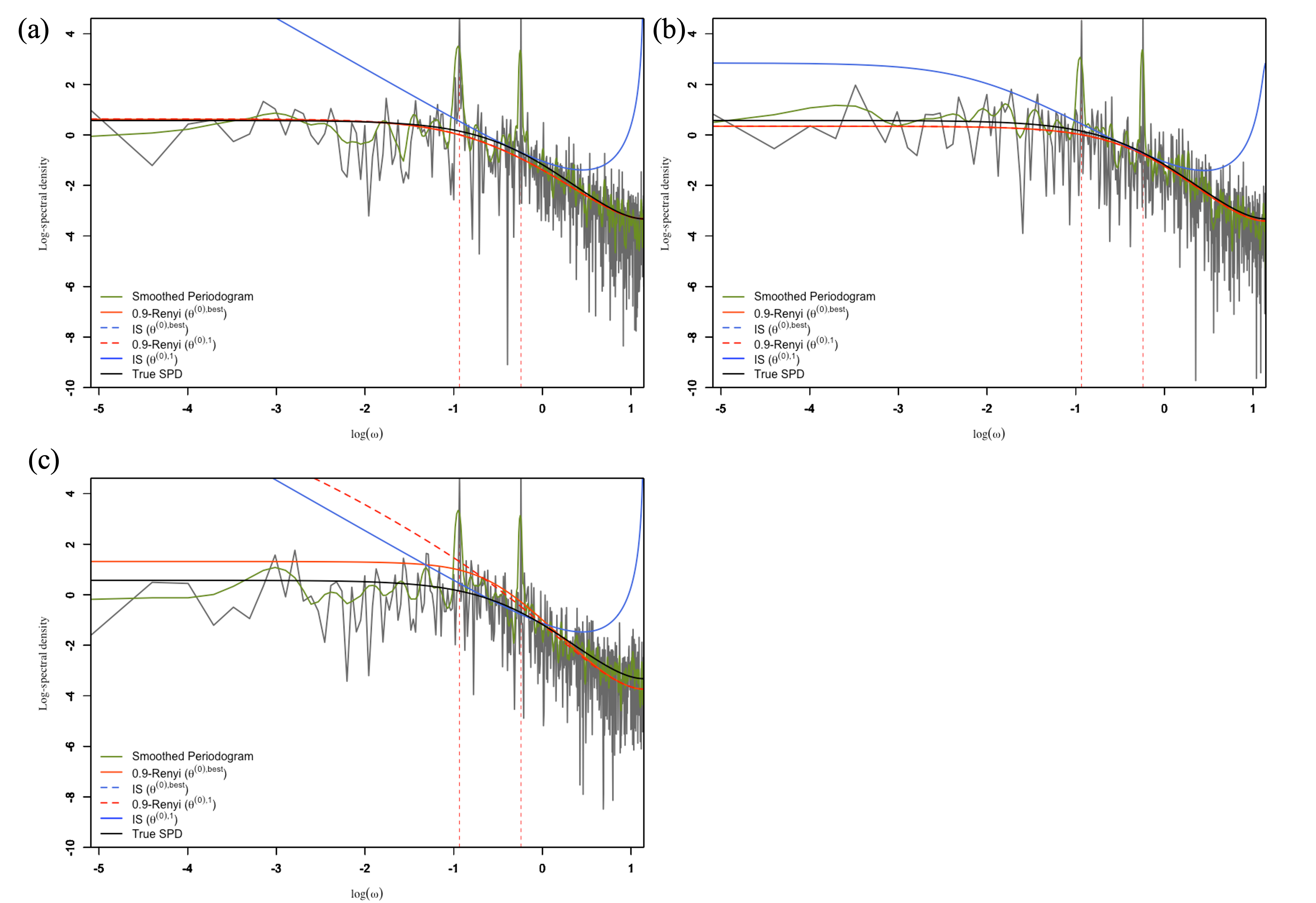}
    \caption{  Spectral 
    %\reviseend 
    densities with the estimates plugged in for the AR(2) model with the trend. 
    Each panel shows the results for different realizations of the periodogram. The legend is the same as in Figure \ref{Fig_SPD1_AR2}.
    The red dot lines represent the points of the outliers. 
    We set $\theta^{(0),1}$ to the initial value.
    }%Trigonometric trends~(outliers). The red dot lines represent the points of the outliers in the frequency domain.}
    \label{Fig_SPD2_AR2}
\end{figure}

\begin{figure}[h]
    \centering
    %\includegraphics[scale=0.43]{Fig_histITER_1.png}
    %\revisebegin
    \includegraphics[scale=0.35]{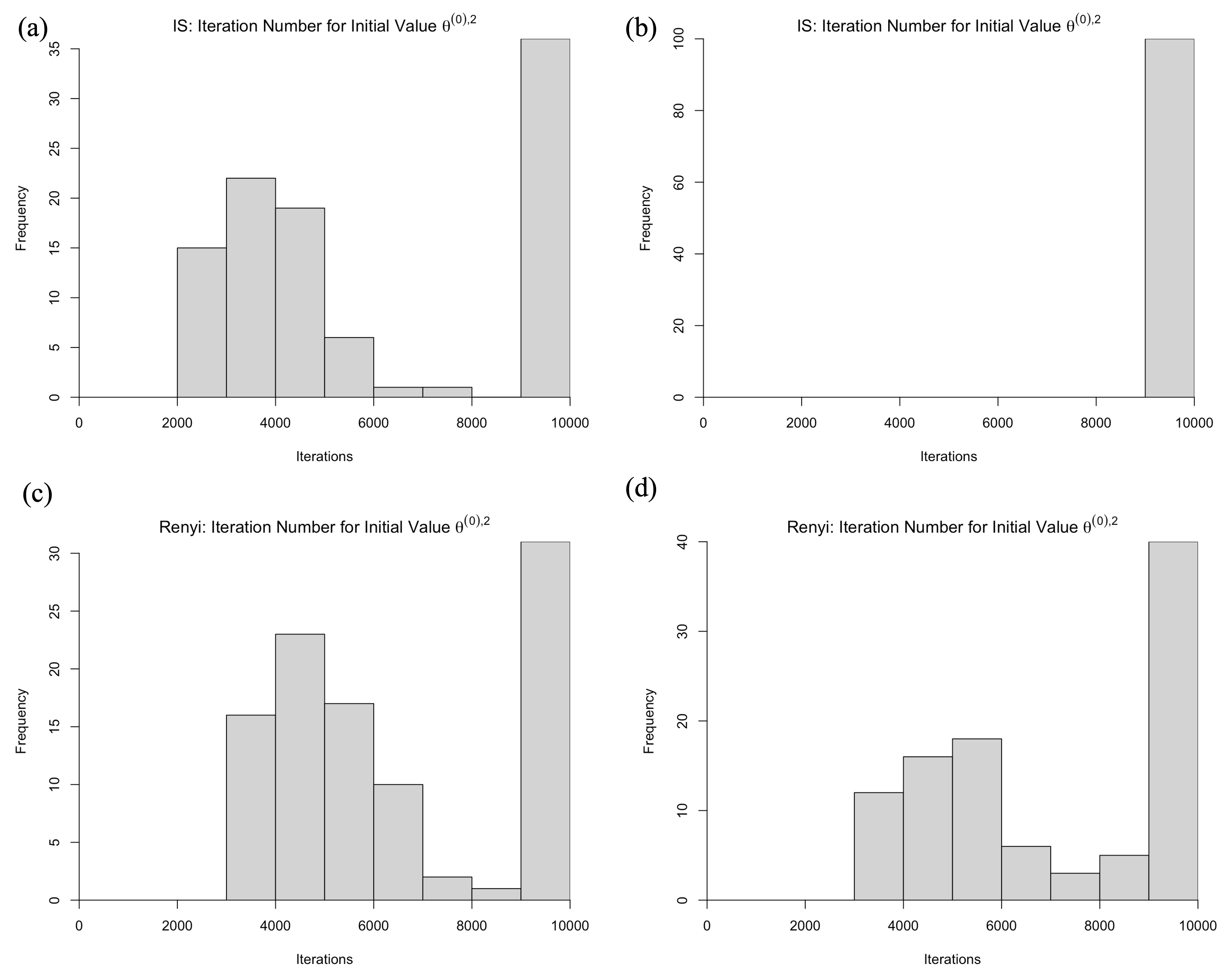}
    \caption{Histograms of the numbers of iterations to converge in the AR(2) model, where the initial value is $\theta^{(0),2}$. %\sout{$\theta^{(0),1}$}. 
    We conducted each simulation 100 times and plotted the histogram. (a) The Itakura--Saito divergence minimization without any trend. (b) The Itakura--Saito divergence minimization with a trend.
    (c) The $\alpha=0.9$-R\'{e}nyi divergence minimization without any trend. (d) The $\alpha=0.9$-R\'{e}nyi divergence minimization with a trend.}\label{Fig_histITER_1}
    %\reviseend
\end{figure}
%\reviseend

Figures~\ref{Fig_SPD1_AR2} and \ref{Fig_SPD2_AR2}
display the estimated spectral densities.
In each figure, 
the solid gray curve represents
the periodogram $I_{n}(\omega)$;
the solid green curve represents the smoothed periodogram $I_{n}^{S}(\omega)$;
the solid black curve represents
the spectral density $S^{\mathrm{AR2}}_{\theta}(\omega)$ 
with the true value $\theta^{*}$ plugged in;
the solid blue curve represents
the spectral density with the minimum Itakura--Saito divergence estimate plugged in;
the solid red curve represents
the spectral density with the minimum spectral R\'{e}nyi divergence ($\alpha=0.9$) estimate plugged in.
The dashed curves represent the estimates starting from the initial values that minimize the objective functions, respectively. Here the periodogram is used 
%\revisebegin
%\sout{instead of the smoothed one,} 
%\reviseend
in these figures.
Figures~\ref{Fig_SPD1_AR2} and \ref{Fig_SPD2_AR2} show that the spectral density $S^{\mathrm{AR2}}_{\widehat{\theta}_{n}}$ with the minimum spectral R\'{e}nyi divergence ($\alpha=0.9$) (resp.~IS divergence) estimate $\widehat{\theta}_{n}$ plugged in fits (resp. does not fit) well to the spectral density $S^{\mathrm{AR2}}_{\theta^{*}}$ 
with the true value $\theta^{*}$ plugged in, which implies that the minimizers of the IS divergence with the periodogram/smoothed periodogram are very sensitive to the outliers in the frequency domain, but the spectral R\'enyi divergence is robust against them. 
These figures also indicate that the Itakura--Saito divergence estimate is highly sensitive to realizations of the periodogram even without any contamination. The grid search for the initial values does not help improve the performance of the Itakura--Saito divergence estimate.
 In contrast, the spectral R\'{e}nyi estimate is generally less sensitive to realizations of the periodogram.

%\revisebegin
Further, Figure \ref{Fig_histITER_1} displays the histograms of the numbers of iterations required for convergence when the initial value is $\theta^{(0),2}$.
For the R\'{e}nyi-divergence minimization, the vast majority of trials converge within a moderate number of iterations, both with and without outliers.
In contrast, for the Itakura--Saito divergence minimization, the convergence behavior deteriorates markedly in the presence of outliers, where essentially all trials reach the maximum number of iterations without convergence.
These findings are consistent with the optimization trajectories illustrated in Figure \ref{Fig_Renyi_demonstrate}.
%\reviseend

\subsection{The Brune spectral model with attenuation}
\label{subsec: Brune}

Next, consider the Brune spectral model \cite{Brune1970} with attenuation \cite{AkiandRichard} for velocity waveforms:
\begin{align*}
S^{\mathrm{BA}}_{\theta}(\omega):=\omega^{2}\left[\frac{\sigma^{2}}{\{1+(\omega/\omega_{c})^{2}\}^{2}}
\right]\exp(-\omega/Q),\;\; \theta=(\sigma,\omega_{c},Q). %\in\Theta.
\end{align*}
The Brune spectral model $\sigma^{2}/\{1+(\omega/\omega_{c})^{2}\}^{2}$ for displacement waveforms  is the spectral density model of Mat\'{e}rn process \cite{Matern1960} with the spectral decay rate fixed to 2
and often used in the earthquake source estimation \cite{CalderoniandAbercrombie2023,Yoshimitsuetal2023}.
The attenuation factor $\exp(-\omega/Q)$ represents anelasticity of the medium through which the seismic wave propagates.
The multiplicative factor $\omega^{2}$ comes from the time derivative that converts a displacement waveform to a velocity waveform.

\begin{figure}[h!]
    \centering
\includegraphics[scale=0.45]{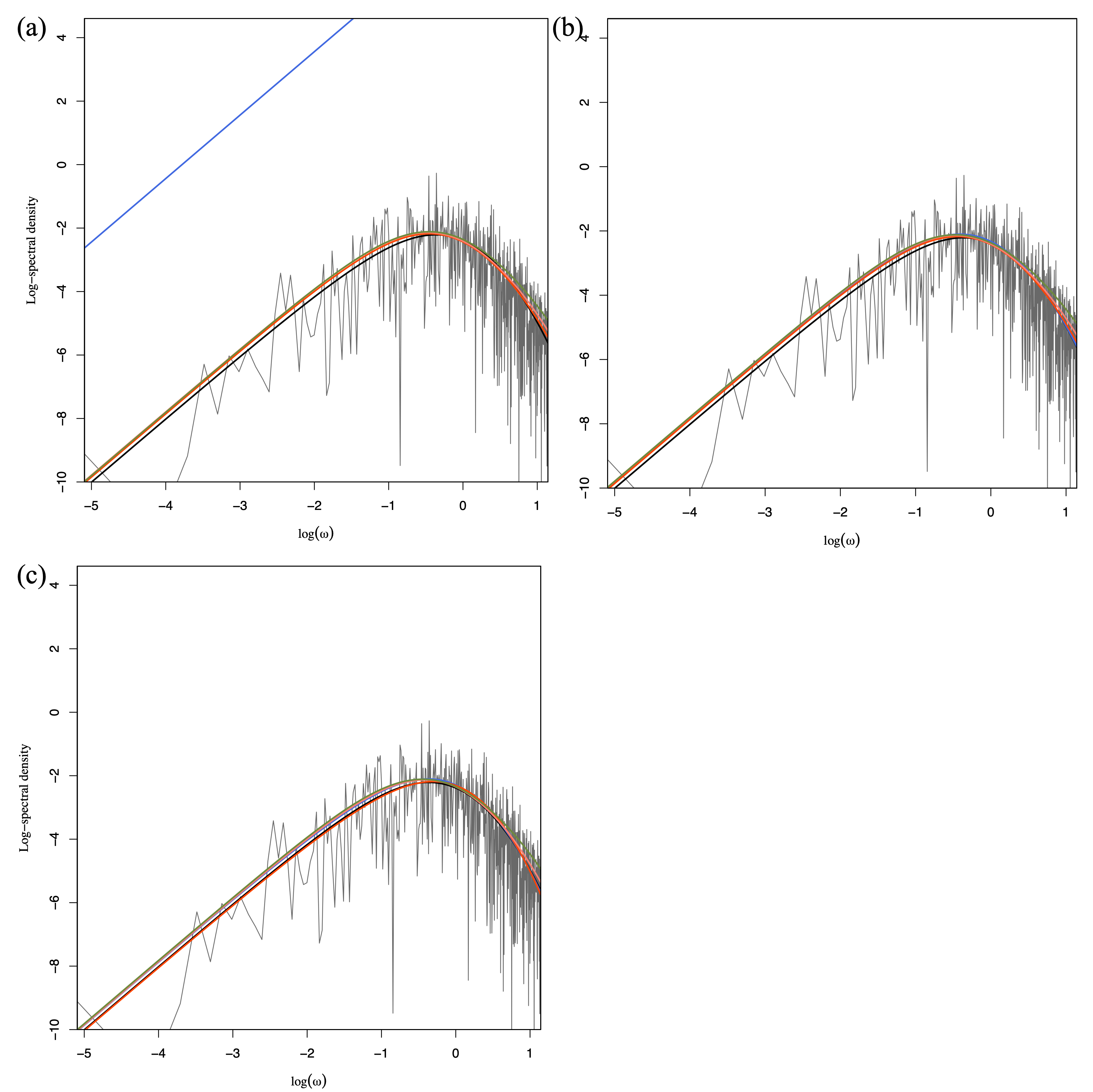}
    \caption{Spectral densities with the estimates plugged in for the Brune model without any trend. The gray curve is the periodogram. The true spectral density is colored in black. Spectral densities based on the spectral R\'{e}nyi divergence ($\alpha=0.5,0.75,0.9$) are colored in red, salmon pink, and green,respectively. The spectral density based on the Itakura--Saito divergence is colored in blue. (a) the result based on the initial value $\theta^{(0),1}$,
    (b) the result based on the initial value $\theta^{(0),2}$,
    (c) the result based on the initial value $\theta^{(0),3}$.
    }
    \label{Fig_SPD1_Brune}
\end{figure}

\begin{figure}[h!]
    \centering
\includegraphics[scale=0.45]{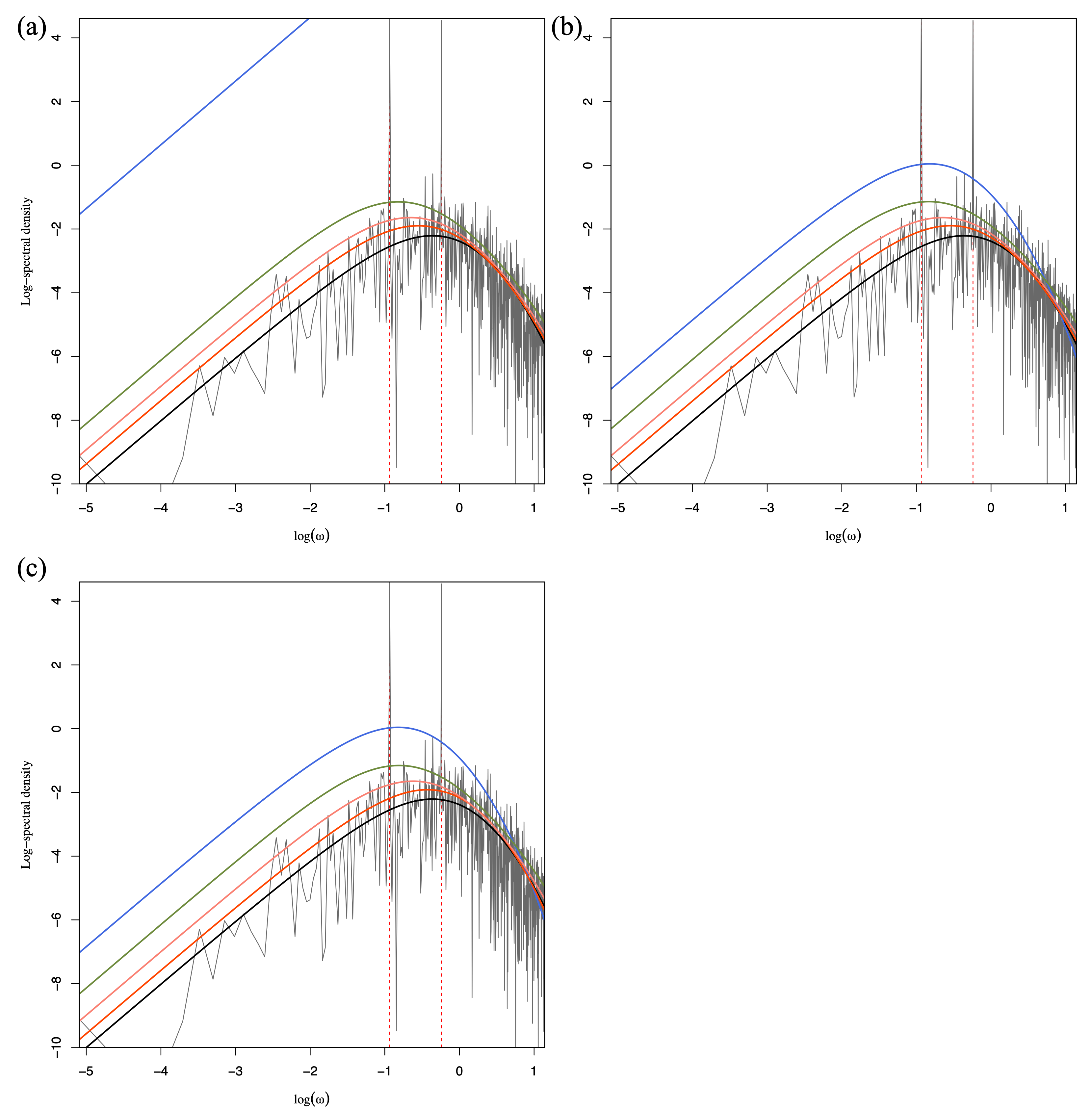}
    \caption{Spectral densities with the estimates plugged in for the Brune model with the trigonometric trends. The red dashed lines denote the frequencies at which the outliers are injected.
    The gray curve is the periodogram. The true spectral density is colored in black. Spectral densities based on the spectral R\'{e}nyi divergence ($\alpha=0.5,0.75,0.9$) are colored in red, salmon pink, and green,respectively. The spectral density based on the Itakura--Saito divergence is colored in blue.
    (a) the result based on the initial value $\theta^{(0),1}$,
    (b) the result based on the initial value $\theta^{(0),2}$,
    (c) the result based on the initial value $\theta^{(0),3}$.
    }
    \label{Fig_SPD2_Brune}
\end{figure}

\begin{table}[h!]
\caption{
Mean values of biases with standard deviations for the Brune model without any trend. Values closest to zero are underlined.
R\'{e}nyi is abbreviated as R; Itakura--Saito is abbreviated as IS; Initial value is abbreviated as Init.
}
\centering
\begin{tabular}{r|crrr}
  & Init& $\hat{\sigma}-\sigma^{*}$ & $\hat{\omega_{c}}-\omega_{c}^{*}$ & $\hat{Q}-Q^{*}$ \\ 
  \hline
R ($\alpha=0.50$) & 
$\theta^{(0),1}$ & -0.03 ($\pm$0.05)  & \underline{-0.12} ($\pm$0.07)& \underline{0.30} ($\pm$0.15) \\ 
  R ($\alpha=0.75$) & $\theta^{(0),1}$  & \underline{0.01} ($\pm$0.06) &  -0.17 ($\pm$0.06) & 0.54 ($\pm$0.22)\\ 
  R ($\alpha=0.90$) &$\theta^{(0),1}$ & \underline{0.01} ($\pm$0.06)& -0.20 ($\pm$0.06) & 0.90 ($\pm$0.33)\\ 
  IS with $I_{n}$  & $\theta^{(0),1}$ &42.6 ($\pm$2.66) & 893 ($\pm$50.5) & 48.2 ($\pm$2.99)\\ 
  IS with $I_{n}^{S}$ & $\theta^{(0),1}$ & 3372 ($\pm$5318) & 67386 ($\pm$106260) & 5259 ($\pm$8327)\\ 
   \hline  \hline
R ($\alpha=0.50$) & 
$\theta^{(0),2}$ & -0.03 ($\pm$0.05)  & -0.10 ($\pm$0.07)& 0.23 ($\pm$0.15) \\ 
  R ($\alpha=0.75$) & $\theta^{(0),2}$  & \underline{0.002} ($\pm$0.06) &  -0.17 ($\pm$0.07) & 0.50 ($\pm$0.21)\\ 
  R ($\alpha=0.90$) &$\theta^{(0),2}$ & 0.01 ($\pm$0.06)& -0.23 ($\pm$0.05) & 1.05 ($\pm$0.34)\\ 
  IS with $I_{n}$  & $\theta^{(0),2}$ & 0.06 ($\pm$0.10) & \underline{-0.04} ($\pm$0.13) & \underline{0.06} ($\pm$0.17)\\ 
  IS with $I_{n}^{S}$ & $\theta^{(0),2}$ & -0.29 ($\pm$0.30) & 1.95 ($\pm$1.89) & 1.86 ($\pm$0.98) \\ 
   \hline  \hline
R ($\alpha=0.50$) & 
$\theta^{(0),3}$ & -0.10 ($\pm$0.04)  &  0.51 ($\pm$0.13)&  -0.30 ($\pm$0.06) \\ 
  R ($\alpha=0.75$) & $\theta^{(0),3}$  & -0.05 ($\pm$0.06) &  \underline{0.08} ($\pm$0.19) &  \underline{0.05} ($\pm$0.20)\\ 
  R ($\alpha=0.90$) &$\theta^{(0),3}$ & \underline{-0.007} ($\pm$0.06)& -0.15 ($\pm$0.13) &  0.72 ($\pm$0.35)\\ 
  IS with $I_{n}$  & $\theta^{(0),3}$ & 0.10 ($\pm$0.07) & 0.42 ($\pm$0.26) &  -0.28 ($\pm$0.13)\\ 
  IS with $I_{n}^{S}$ & $\theta^{(0),3}$ & -0.29 ($\pm$0.15) & 2.09 ($\pm$1.68) &  1.50 ($\pm$0.72)\\ 
   \hline 
\end{tabular}
\label{tab: Brune without trend}
\end{table}

\begin{table}[h!]
\caption{
Mean values of biases with standard deviations for the Brune model with the trigonometric trends. Values closest to zero are underlined.
R\'{e}nyi is abbreviated as R; Itakura--Saito is abbreviated as IS; Initial value is abbreviated as Init.
}
\centering
\begin{tabular}{r|crrr}
  & Init & $\hat{\sigma}-\sigma^{*}$ & $\hat{\omega_{c}}-\omega_{c}^{*}$ & $\hat{Q}-Q^{*}$ \\ 
  \hline
R ($\alpha=0.50$) & 
$\theta^{(0),1}$ & \underline{0.22} ($\pm$0.07)  & \underline{-0.22} ($\pm$0.05)& \underline{0.27} ($\pm$0.15) \\ 
  R ($\alpha=0.75$) & $\theta^{(0),1}$  & 0.53 ($\pm$0.08) & -0.34 ($\pm$0.04) & 0.47 ($\pm$0.19)\\ 
  R ($\alpha=0.90$) &$\theta^{(0),1}$ & 1.45 ($\pm$0.10) & -0.48 ($\pm$0.03) &  0.74 ($\pm$0.31)\\ 
  IS with $I_{n}$  & $\theta^{(0),1}$ & 75.8 ($\pm$2.85) &  1524 ($\pm$55.86) &  60.0 ($\pm$3.04)\\ 
 IS with $I_{n}^{S}$ & $\theta^{(0),1}$ &  3405 ($\pm$5318) &  68019 ($\pm$106263) &  5271 ($\pm$8327)\\ 
   \hline  \hline
R ($\alpha=0.50$) & 
$\theta^{(0),2}$ & \underline{0.21} ($\pm$0.07)  &   \underline{-0.20} ($\pm$0.06)& \underline{0.23} ($\pm$0.14) \\ 
  R ($\alpha=0.75$) & $\theta^{(0),2}$  & 0.53 ($\pm$0.08) &   -0.33 ($\pm$0.04) & 0.45 ($\pm$0.19)\\ 
  R ($\alpha=0.90$) &$\theta^{(0),2}$ &  1.48 ($\pm$0.10)&  -0.49 ($\pm$0.03) & 0.79 ($\pm$0.33)\\ 
  IS with $I_{n}$  & $\theta^{(0),2}$ & 3.92 ($\pm$0.07) &  -0.35 ($\pm$0.02) & -0.40 ($\pm$0.02)\\ 
  IS with $I_{n}^{S}$ & $\theta^{(0),2}$ & 2.83 ($\pm$0.97) &  -0.09 ($\pm$0.97) &  2.30 ($\pm$0.98) \\ 
   \hline  \hline
R ($\alpha=0.50$) & 
$\theta^{(0),3}$ & \underline{0.08} ($\pm$0.06)  &  \underline{0.24} ($\pm$0.15)& \underline{-0.25} ($\pm$0.10) \\ 
  R ($\alpha=0.75$) & $\theta^{(0),3}$  & 0.46 ($\pm$0.08) &   -0.27 ($\pm$0.07) & 0.26 ($\pm$0.16)\\ 
  R ($\alpha=0.90$) &$\theta^{(0),3}$ & 1.40 ($\pm$0.09)&  -0.47 ($\pm$0.03) & 0.69 ($\pm$0.29) \\ 
  IS with $I_{n}$  & $\theta^{(0),3}$ & 3.90 ($\pm$0.07) &  -0.34 ($\pm$0.02) & -0.40 ($\pm$0.02)\\ 
  IS with $I_{n}^{S}$ & $\theta^{(0),3}$ & 2.80 ($\pm$0.95) & -0.13 ($\pm$0.59) & 2.01 ($\pm$0.86)\\ 
   \hline 
\end{tabular}
\label{tab: Brune with trend}
\end{table}

For the experiments,
we use the true value 
$\theta^{\ast}:=(\sigma^{\ast},\omega_{c}^{\ast},Q^{\ast})=(1,1,1)$.
We generate a time series $\mathbf{x}_{n=2^{10}}$ with its spectral density $S_{\theta^{*}}^{\mathrm{BA}}$
100 times,
and if we consider outliers in the frequency domain,
we add two trigonometric trends to the time series as
\[
x_{t}^{o}=x_{t}+\sqrt{8\pi z_{1}/n}\sin(A_{1}t)+\sqrt{8\pi z_{2}/n}\sin(A_{2}t),\,\,
t=1,\ldots,n
\]
with $z_{1}=z_{2}=100$, $A_{1}=\pi/4$ and $A_{2}=\pi/8$.
We report the results based on different values of $z_{1}=z_{2}$ in Appendix \ref{appendix: additional examples}.
We use three different initial values:
$\theta^{(0),1}=(1,0.1,1)$,
$\theta^{(0),2}=(1,1,1)$,
and
$\theta^{(0),3}=(1,2,1)$.
Setting $\sigma^{*}$ and $Q^{*}$ as initial values does not impact on the behavior of the minimum spectral R\'{e}nyi divergence estimates; we report the results based on different initial values 
$\theta^{(0),4}=(0.1,0.1,0.1)$
and 
$\theta^{(0),5}=(2,2,2)$ in Appendix \ref{appendix: additional examples}.

Figures~\ref{Fig_SPD1_Brune} and \ref{Fig_SPD2_Brune} display the estimated spectral densities. In each figure, the following representations are used:
\begin{itemize}
\item The solid gray curve represents the periodogram $I_{n}$;
\item the solid black curve denotes the true spectral density $S^{\mathrm{BA}}_{\theta^{*}}$;
\item the solid blue curve illustrates the spectral density with the minimum Itakura--Saito divergence estimate plugged in;
\item the solid red, salmon pink, and green curves display the spectral densities with the minimum spectral R\'{e}nyi divergence ($\alpha=0.9,0.75,0.5$) estimates plugged in, respectively.
\end{itemize}

For any initial value and regardless of the existence of the trend, the spectral density based on the spectral $\alpha=0.5$-R\'{e}nyi divergence aligns well with the true spectral density. The minimum spectral R\'{e}nyi divergence estimate with any $\alpha\in\{0.5,0.75,0.9\}$ remains stable with respect to the choice of the initial value. Conversely, the minimum Itakura--Saito divergence estimate shows sensitivity to this choice.

Tables \ref{tab: Brune without trend} and \ref{tab: Brune with trend} summarize the estimation results. 
Comparing the results based on different initial values, we find that the minimum Itakura--Saito divergence estimates are sensitive to the initial value of the optimization regardless of the existence of the smoother,
while the spectral R\'enyi divergence yields estimation results stable to the choice of the initial value. 
The minimum Itakura--Saito divergence estimate with the smoothed periodogram performs much worse than that with the periodogram. 
In the presence of outliers in the frequency domain, 
the minimum R\'{e}nyi divergence estimate with $\alpha=0.5$ performs the best regardless of the choice of the initial value in this example.

\section{Discussions}

Finally, we discuss three important directions for future research.

The first point is the choice of $\alpha$.
The value of $\alpha$ controls the trade-off between the efficiency without the presence of outliers and the robustness even for larger trends.
For an objective selection of $\alpha$,
constructing an information criterion is one possible direction, but requires sufficient discussions and experiments.
In practice, 
monitoring the behaviors of the minimum spectral R\'{e}nyi divergence estimates for several values of $\alpha$ (say, $\alpha=0.5,0.75,0.9$) would be a workaround.

The second point is the choice of optimization schemes.
In this paper, we employed a simple gradient descent scheme only as a baseline illustrative method, and we do not claim that naive gradient descent is the most appropriate choice in practice. As is well known, the iteration number, tolerance, and step-size selection in naive gradient descent can be highly sensitive and not objectively determined.
There are at least two ways to avoid these issues:
\begin{enumerate}
    \item \textbf{Use of line-search methods.} 
    One can incorporate standard line-search strategies such as Wolfe or Armijo conditions (see Proposition~2 and Remark~3), or alternatively employ recent line-search-free first-order methods \cite{malitsky2020adaptive, yagishita2025simple}, which often provide more stable convergence behavior without requiring excessively large iteration limits.
    \item \textbf{Use of quasi-Newton methods.}  
    Quasi-Newton algorithms such as 
    the Broyden--Fletcher--Goldfarb--Shanno (BFGS) algorithm or L-BFGS \cite{liu1989limited} are typically more efficient in terms of optimization and are widely preferred by practitioners. Indeed, in a large-scale application of spectral R\'{e}nyi divergence minimization to real geophysical data \cite{kano2025spatiotemporal}, we adopted the L-BFGS method precisely to avoid the practical difficulties associated with naive first-order methods such as gradient descent.
\end{enumerate}
Constructing a unified theory that accommodates line-search-free first-order methods and quasi-Newton methods 
would be an important research direction.

The third one is the impact of the choice of an initial value.
Numerical experiments presented in Section \ref{sec:numerical} suggest that the minimum spectral R\'{e}nyi divergence estimates are robust to the choice of an initial value, while the minimum Itakura--Saito divergence estimate is sensitive to it.
These robustness and sensitivity are confirmed in the other examples such as autoregressive models.
The structure of the gradient $\mathcal{G}_{\alpha}(\theta\,;\,\tilde{I}_{n})$
in the comparison with that of $\mathcal{G}_{\mathrm{IS}}(\theta\,;\,\tilde{I}_{n})$ 
can give a clue to this, and so detailed analysis of these gradient landscapes would be one of the interesting research directions.

\section{Acknowledgements}
The authors would like to thank the Associate Editor and the referee for their constructive comments, which improved the quality of the paper.
The authors would like to thank Akifumi Okuno, Mirai Tanaka, Kei Kobayashi, Yuta Koike, Tomoyuki Higuchi, Hironori Fujisawa, Junho Yang, Masayuki Kano, and Shotaro Yagishita for their comments.
This work is supported by JSPS KAKENHI (19K20222, 21H05205, 21K12067, 23K11024), MEXT (JPJ010217), and ``Strategic Research Projects'' grant (2022-SRP-13) from ROIS (Research Organization of Information and Systems). 
The simulation code is available at \url{https://github.com/t-tetsuya/Spectral_Renyi.git}.

%\clearpage
%\newpage
\appendix
%{\Large Supplement to ``On robustness of spectral R\'{e}nyi divergence" by Tetsuya Takabatake and Keisuke Yano}
\section*{Supplement to ``On robustness of spectral R\'{e}nyi divergence" by Tetsuya Takabatake and Keisuke Yano}

This supplement presents
\begin{itemize}
    \item Proof of Theorem \ref{thm: Gamma meets Renyi} in Appendix \ref{appendix: proof of Gamma};
    \item Proof of Theorem \ref{thm: success Renyi} in Appendix \ref{appendix: proof of thm};
    \item Proof of Proposition \ref{prop: stability of line search} in Appendix \ref{appendix: Proof of line search};
    \item Proof of Proposition \ref{prop: failure Itakura Saito} in Appendix \ref{appendix: proof of Prop}; 
    \item %\revisebegin 
    Additional comparison of optimization paths in Appendix \ref{appendix: additional comparison}; and %\reviseend
    \item Additional simulation studies in Appendix \ref{appendix: additional examples}.
\end{itemize}

\section{Proof of Theorem \ref{thm: Gamma meets Renyi}}
\label{appendix: proof of Gamma}

We first introduce the notation used in the proof. 
For each spectral density function $S$ of a stationary time series %\sout{$\{X_{t}\}_{t\in\mathbb{Z}}$,}
%\revisebegin
$(X_{t})_{t\in\mathbb{Z}}$, 
%\reviseend
we write $\Sigma_{n}(S)$ as the $n\times n$-Toeplitz matrix whose $(s,t)$th element is given by
\begin{align*}
    R(t-s)=\int_{-\pi}^{\pi}e^{\sqrt{-1}(t-s)\omega}S(\omega)\,\mathrm{d}\omega,\ \ s,t=1,\ldots,n,
\end{align*}
which corresponds to the covariance matrix of $(X_{1},X_{2},\cdots,X_{n})^{\top}$. 
Then, the probability density function of $(X_{1},\cdots,X_{n})^{\top}$ is expressed by
\begin{align}\label{def:pdf_mGauss}
    p_{S}(\mathbf{x})=\frac{1}{\sqrt{(2\pi)^{n}\det[\Sigma_{n}(S)]}}\exp\left(-\frac{1}{2}\mathbf{x}^{\top}\Sigma_{n}(S)^{-1}\mathbf{x}\right),\ \ \mathbf{x}\in\mathbb{R}^{n}.
\end{align}
Let $\gamma\in(0,\infty)$. 
First introduce the $\gamma$-cross entropy %\sout{$d_{\gamma}(p_{S},p_{\widetilde{S}})$}
%\revisebegin
$d_{\gamma}[\,p_{S}\,:\,p_{\widetilde{S}}\,]$ 
%\reviseend
between two probability density functions $p_{S}$ and $p_{\widetilde{S}}$, which are not necessarily Gaussian probability density functions, defined by %\sout{$d_{\gamma}[p_{S}:p_{\widetilde{S}}]$}
\begin{align*}
    %\revisebegin
    d_{\gamma}[\,p_{S}\,:\,p_{\widetilde{S}}\,]
    %\reviseend
    :=-\frac{1}{\gamma}\log\int_{\mathbb{R}^{n}}p_{S}(\mathbf{x})p_{\widetilde{S}}(\mathbf{x})^{\gamma}\,\mathrm{d}\mathbf{x}
    +\frac{1}{1+\gamma}
    \log\int_{\mathbb{R}^{n}}p_{\widetilde{S}}(\mathbf{x})^{1+\gamma}\intd\mathbf{x}.
\end{align*}
Then the $\gamma$-divergence can be written as %\sout{$D_{\gamma}[p_{S}:p_{\widetilde{S}}]=d_{\gamma}[p_{S}~:p_{\widetilde{S}}]-d_{\gamma}[p_{S}:p_{S}]$}
%\revisebegin
$G_{\gamma}[\,p_{S}\,:\,p_{\widetilde{S}}\,]=d_{\gamma}[\,p_{S}\,:\,p_{\widetilde{S}}\,]-d_{\gamma}[\,p_{S}\,:\,p_{S}\,]$.
%\reviseend 
\\

We now return to the proof. We first compute the $\gamma$-cross entropy %\sout{$d_{\gamma}[p_{S}:p_{\widetilde{S}}]$}
%\revisebegin
$d_{\gamma}[\,p_{S}\,:\,p_{\widetilde{S}}\,]$
%\reviseend
for Gaussian probability density functions $p_{S}$ and $p_{\widetilde{S}}$, and then the conclusion directly follows from Szeg\"o's limit theorem \cite{Szego_1920,Kolmogorov_1941}. 
Since $S$ and $\tilde{S}$ are positive on $(0,\pi)$, we can show that $\Sigma_{n}(S)$ and $\Sigma_{n}(\widetilde{S})$ are positive definite so that $\Sigma_{n}(\widetilde{S})+\gamma \Sigma_{n}(S)$ is also positive definite for any $\gamma>0$, which implies that $\Sigma_{n}(\widetilde{S})+\gamma \Sigma_{n}(S)$ is invertible. 
Then we obtain
\begin{align*}
    \int_{\mathbb{R}^{n}}p_{S}(\mathbf{x})p_{\widetilde{S}}(\mathbf{x})^{\gamma}\intd\mathbf{x} 
    &=\left\{\pi^{n(1+\gamma)}\det[\Sigma_{n}(S)]\det[\Sigma_{n}(\widetilde{S})]^{\gamma}\right\}^{-\frac{1}{2}}\\
    &\qquad\quad
    \int_{\mathbb{R}^{n}}\exp\left(-\frac{1}{2}\mathbf{x}^{\top}[\Sigma_{n}(S)^{-1}+\gamma\Sigma_{n}(\widetilde{S})^{-1}]\mathbf{x}\right)\intd\mathbf{x} \\
    &=\left\{\pi^{n(1+\gamma)}\det[\Sigma_{n}(S)]\det[\Sigma_{n}(\widetilde{S})]^{\gamma}\right\}^{-\frac{1}{2}}\\
    &\qquad\quad
    \int_{\mathbb{R}^{n}}\exp\left(-\frac{1}{2}\mathbf{x}^{\top}\Sigma_{n}(S)^{-1}[\Sigma_{n}(\widetilde{S})+\gamma \Sigma_{n}(S)]\Sigma_{n}(\widetilde{S})^{-1}\mathbf{x}\right)\intd\mathbf{x} \\
    &=\left\{\frac{\pi^{n}\det[\Sigma_{n}(\widetilde{S})(\Sigma_{n}(\widetilde{S})+\gamma \Sigma_{n}(S))^{-1}\Sigma_{n}(S)]}{\pi^{n(1+\gamma)}\det[\Sigma_{n}(S)]\det[\Sigma_{n}(\widetilde{S})]^{\gamma}}\right\}^{\frac{1}{2}}  \\
    &=\left\{\pi^{-n\gamma}\det[\Sigma_{n}(\widetilde{S})]^{1-\gamma}
    \det[\Sigma_{n}(\widetilde{S})+\gamma \Sigma_{n}(S)]^{-1}\right\}^{\frac{1}{2}}.
\end{align*}
Similarly, we also obtain
\begin{align*}
    \int_{\mathbb{R}^{n}}p_{\widetilde{S}}(\mathbf{x})^{1+\gamma}\intd\mathbf{x}
    &=\left\{\pi^{-n\gamma}\det[\Sigma_{n}(\widetilde{S})]^{1-\gamma}
    \det[(1+\gamma)\Sigma_{n}(\widetilde{S})]^{-1}\right\}^{\frac{1}{2}}\\
    &=\left\{\pi^{-n\gamma}(1+\gamma)^{-n}\det[\Sigma_{n}(\widetilde{S})]^{-\gamma}\right\}^{\frac{1}{2}}.
\end{align*}
Combining these yields %\sout{$2d_{\gamma}[p_{S}:p_{\widetilde{S}}]$}
\begin{align*}
    %\revisebegin
    2d_{\gamma}[\,p_{S}\,:\,p_{\widetilde{S}}\,]
    %\reviseend
    =&-\frac{1}{\gamma}\log\left[\pi^{-n\gamma}\det[\Sigma_{n}(\widetilde{S})]^{1-\gamma}
    \det[\Sigma_{n}(\widetilde{S})+\gamma \Sigma_{n}(S)]^{-1}\right] \\
    &\qquad+\frac{1}{1+\gamma}%\prod_{j=1}^{n}
    \log\left[\pi^{-n\gamma}(1+\gamma)^{-n}\det[\Sigma_{n}(\widetilde{S})]^{-\gamma}\right] \\
    =&n\log(\pi) -\frac{n}{1+\gamma}\log[(1+\gamma)\pi^{\gamma}] \\
    &\qquad+\left[
    \frac{1}{\gamma}\log\det[\Sigma_{n}(\widetilde{S})+\gamma \Sigma_{n}(S)]
    -\frac{1}{\gamma(1+\gamma)}\log\det[\Sigma_{n}(\widetilde{S})]
    \right].
\end{align*}
Putting $S$ into $\widetilde{S}$, we obtain %\sout{$2d_{\gamma}[p_{S}:p_{S}]$}
\begin{align*}
    %\revisebegin
    2d_{\gamma}[\,p_{S}\,:\,p_{S}\,]
    %\reviseend
    &=n\log(\pi) -\frac{n}{1+\gamma}\log[(1+\gamma)\pi^{\gamma}] 
    +\frac{1}{\gamma}\log\det[(1+\gamma)\Sigma_{n}(S)] \\
    &\qquad\qquad\qquad\qquad\qquad
    -\frac{1}{\gamma(1+\gamma)}\log\det[\Sigma_{n}(S)] \\
    &=n\log(\pi) -\frac{n}{1+\gamma}\log[(1+\gamma)\pi^{\gamma}] 
    +\frac{n}{\gamma}\log(1+\gamma)
    +\frac{1}{1+\gamma}\log\det[\Sigma_{n}(S)].
\end{align*}
Thus we get %\sout{$2D_{\gamma}[p_{S}:p_{\widetilde{S}}]$}
\begin{align*}
    %\revisebegin
    2G_{\gamma}[\,p_{S}\,:\,p_{\widetilde{S}}\,]
    %\reviseend
    &=\frac{1}{\gamma}\Bigg{[}
    \log\det[\Sigma_{n}(\widetilde{S})+\gamma \Sigma_{n}(S)]
    -\frac{1}{1+\gamma}\log\det[\Sigma_{n}(\widetilde{S})]\\
    &\qquad\qquad\qquad\qquad\qquad\qquad
    -\log(1+\gamma)^{n}
    -\frac{\gamma}{1+\gamma}\log\det[\Sigma_{n}(S)]
    \Bigg{]} \\
    &=\frac{1}{\gamma}\Bigg{[}
    \log\det\left[\frac{\gamma}{1+\gamma} \Sigma_{n}(S) +\frac{1}{1+\gamma}\Sigma_{n}(\widetilde{S})\right]\\
    &\qquad\qquad\qquad\qquad
    -\frac{\gamma}{1+\gamma}\log\det[\Sigma_{n}(S)]
    -\frac{1}{1+\gamma}\log\det[\Sigma_{n}(\widetilde{S})]
    \Bigg{]} \\
    &=\frac{1}{\gamma}\left[
    \log\det[\Sigma_{n}(\bar{S}^{(\gamma)})] -\frac{\gamma}{1+\gamma}\log\det[\Sigma_{n}(S)]
    -\frac{1}{1+\gamma}\log\det[\Sigma_{n}(\widetilde{S})]
    \right],
\end{align*}
where $\bar{S}^{(\gamma)}(\omega):=\frac{\gamma}{1+\gamma} S(\omega) +\frac{1}{1+\gamma}\widetilde{S}(\omega)$. 
Using Szeg\"o's limit theorem \cite{Szego_1920,Kolmogorov_1941}, we conclude
%\sout{$D_{\gamma}[p_{S}:p_{\widetilde{S}}]$, $D_{\gamma}[S:\widetilde{S}]$}
\begin{align*}
    \lim_{n\to\infty}\frac{2}{n}
    %\revisebegin
    G_{\gamma}[\,p_{S}\,:\,p_{\widetilde{S}}\,]
    %\reviseend
    &=\frac{1}{2\pi\gamma}\int_{-\pi}^{\pi}\left[
    \log{\bar{S}^{(\gamma)}(\omega)} -\frac{\gamma}{1+\gamma}\log{S(\omega)}
    -\frac{1}{1+\gamma}\log{\widetilde{S}(\omega)} \right]\intd\omega \\
    &= 
    %\revisebegin
    D_{\alpha}[\,S\,:\,\widetilde{S}\,],
    %\reviseend
\end{align*}
%\revisebegin
where $\alpha:=(1+\gamma)^{-1}$. 
%\reviseend
This completes the proof. 
%\qed

\section{Proof of Theorem \ref{thm: success Renyi}}
\label{appendix: proof of thm}

%\begin{proof}
We prove Theorem \ref{thm: success Renyi} by the mathematical induction.

\textbf{Step 1}: Consider $m=1$. Observe
\begin{align*}
&\hat{\theta}^{(1)}_{\alpha}[\tilde{I}^{z,\omega^{*}}_{n}]-
\hat{\theta}^{(1)}_{\alpha}[\tilde{I}_{n}]\\
&=\left\{
\hat{\theta}^{(1)}_{\alpha}[\tilde{I}^{z,\omega^{*}}_{n}]-\theta^{(0)}
\right\}
-\left\{
\hat{\theta}^{(1)}_{\alpha}[\tilde{I}_{n}]
-\theta^{(0)}
\right\}\\
&=
\frac{\gamma_{1} \alpha}{(1-\alpha)n}
\Bigg{(}
\frac{S_{\theta^{(0)}}(\omega^{*})}{\alpha S_{\theta^{(0)}}(\omega^{*})+ (1-\alpha)\tilde{I}_{n}(\omega^{*})}\\
&\qquad\qquad\qquad\qquad-
\frac{S_{\theta^{(0)}}(\omega^{*})}{\alpha S_{\theta^{(0)}}(\omega^{*})+ (1-\alpha)(\tilde{I}_{n}(\omega^{*})+z)}
\Bigg{)}\nabla_{\theta^{(0)}}\log S_{\theta}(\omega^{*})
\\
&=\frac{\gamma_{1}}{(1-\alpha)n}
\frac{\alpha S_{\theta^{(0)}}(\omega^{*})}{\alpha S_{\theta^{(0)}}(\omega^{*})+(1-\alpha)\tilde{I}_{n}(\omega^{*})}
\frac{(1-\alpha)z \nabla_{\theta^{(0)}}S_{\theta}(\omega^{*})}{(1-\alpha)z+(1-\alpha)\tilde{I}_{n}(\omega^{*})+\alpha S_{\theta^{(0)}}(\omega^{*})}.
\end{align*}
This yields
\begin{align}
\left\|\hat{\theta}^{(1)}_{\alpha}[\tilde{I}^{z,\omega^{*}}_{n}]-
\hat{\theta}^{(1)}_{\alpha}[\tilde{I}_{n}]
\right\|
\le 
\frac{1}{n}\frac{\gamma_{1}}{1-\alpha}\|\nabla_{\theta^{(0)}}S_{\theta}(\omega^{*})\|,
\label{eq: comparison of first step}
\end{align}
which proves the assertion for $m=1$.

\textbf{Step 2}: 
Let $m\in\mathbb{N}>1$.
Assume that for any $\varepsilon>0$,
there exists $N\in\mathbb{N}$ such that for $n\ge N$, we have
\begin{align*}
\|\hat{\theta}^{(m)}_{\alpha}[\tilde{I}_{n}^{z,\omega^{*}}]
-
\hat{\theta}^{(m)}_{\alpha}[\tilde{I}_{n}]\|\le \varepsilon.
\end{align*}
By the cancelling technique and by the triangle inequality, we get 
\begin{align*}
\left\|\hat{\theta}^{(m+1)}_{\alpha}[\tilde{I}_{n}^{z,\omega^{*}}]
-
\hat{\theta}^{(m+1)}_{\alpha}[\tilde{I}_{n}]\right\|
&=
\Bigg{\|}\left\{\hat{\theta}^{(m+1)}_{\alpha}[\tilde{I}_{n}^{z,\omega^{*}}]
-
\hat{\theta}^{(m)}_{\alpha}[\tilde{I}_{n}^{z,\omega^{*}}]
\right\} \\
&\qquad%\qquad\qquad\qquad\qquad\qquad\qquad\qquad\qquad
-
\left\{\hat{\theta}^{(m+1)}_{\alpha}[\tilde{I}_{n}]
-
\hat{\theta}^{(m)}_{\alpha}[\tilde{I}_{n}]
\right\}
-\left\{ \hat{\theta}^{(m)}_{\alpha}[\tilde{I}_{n}^{z,\omega^{*}}]
-
\hat{\theta}^{(m)}_{\alpha}[\tilde{I}_{n}]
\right\}
\Bigg{\|}\\
&\le
\gamma_{m}\left\|
\mathcal{G}_{\alpha}\left(\hat{\theta}^{(m)}_{\alpha}[\tilde{I}^{z,\omega^{*}}_{n}]\,;\,
\tilde{I}^{z,\omega^{*}}_{n}
\right)
-
\mathcal{G}_{\alpha}\left(\hat{\theta}^{(m)}_{\alpha}[\tilde{I}_{n}]\,;\,
\tilde{I}_{n}
\right)
\right\|
+\varepsilon.
\end{align*}
Observe that by applying the Taylor theorem to $\nabla_{\theta}\log S_{\theta}(\omega)$
and letting $t$ be some point on some line connecting $\hat{\theta}^{(m)}_{\alpha}[\tilde{I}^{z,\omega^{*}}_{n}]$
to $\hat{\theta}^{(m)}_{\alpha}[\tilde{I}_{n}]$, we get
\begin{align*}
&\mathcal{G}_{\alpha}\left(\hat{\theta}^{(m)}_{\alpha}[\tilde{I}^{z,\omega^{*}}_{n}]\,;\,
\tilde{I}^{z,\omega^{*}}_{n}
\right)
\\
&=
\frac{\alpha}{(1-\alpha)n}
\sum_{\omega\in\Omega_{n}}
\left\{
1-\frac{S_{\hat{\theta}^{(m)}_{\alpha}[\tilde{I}^{z,\omega^{*}}_{n}]}(\omega) }
{\alpha S_{\hat{\theta}^{(m)}_{\alpha}[\tilde{I}^{z,\omega^{*}}_{n}]}(\omega)
+
(1-\alpha)\tilde{I}^{z,\omega^{*}}_{n}(\omega)
}
\right\}\\
&\qquad\qquad\qquad\qquad\qquad\qquad\qquad\qquad\qquad\qquad\qquad\qquad
\nabla_{\hat{\theta}^{(m)}_{\alpha}[\tilde{I}^{z,\omega^{*}}_{n}]}\log S_{\theta}(\omega)\\
&=
\frac{\alpha}{(1-\alpha)n}
\sum_{\omega\in\Omega_{n}}
\left\{
1-\frac{S_{\hat{\theta}^{(m)}_{\alpha}[\tilde{I}^{z,\omega^{*}}_{n}]}(\omega) }
{\alpha S_{\hat{\theta}^{(m)}_{\alpha}[\tilde{I}^{z,\omega^{*}}_{n}]}(\omega)
+
(1-\alpha)\tilde{I}^{z,\omega^{*}}_{n}(\omega)
}
\right\}\nabla_{\hat{\theta}^{(m)}_{\alpha}[\tilde{I}_{n}]}\log S_{\theta}(\omega)\\
&\quad
+\frac{\alpha}{(1-\alpha)n}
\sum_{\omega\in\Omega_{n}}
\left\{
1-\frac{S_{\hat{\theta}^{(m)}_{\alpha}[\tilde{I}^{z,\omega^{*}}_{n}]}(\omega) }
{\alpha S_{\hat{\theta}^{(m)}_{\alpha}[\tilde{I}^{z,\omega^{*}}_{n}]}(\omega)
+
(1-\alpha)\tilde{I}^{z,\omega^{*}}_{n}(\omega)
}
\right\}\nabla_{t}^{2}\log S_{\theta}(\omega)\\
&\qquad\qquad\qquad\qquad\qquad\qquad\qquad\qquad\qquad
(\hat{\theta}^{(m)}_{\alpha}[\tilde{I}^{z,\omega^{*}}_{n}] -
\hat{\theta}^{(m)}_{\alpha}[\tilde{I}_{n}]
)\\
&=
\frac{\alpha}{(1-\alpha)n}
\sum_{\omega\in\Omega_{n}}
\left\{
1-\frac{S_{\hat{\theta}^{(m)}_{\alpha}[\tilde{I}^{z,\omega^{*}}_{n}]}(\omega) }
{\alpha S_{\hat{\theta}^{(m)}_{\alpha}[\tilde{I}^{z,\omega^{*}}_{n}]}(\omega)
+
(1-\alpha)\tilde{I}^{z,\omega^{*}}_{n}(\omega)
}
\right\}\nabla_{\hat{\theta}^{(m)}_{\alpha}[\tilde{I}_{n}]}\log S_{\theta}(\omega)\\
&\quad
+O\left(\frac{2\varepsilon }{(1-\alpha)n}\sup_{\theta\in
%\mbox{\sout{$\Theta$}}
%\revisebegin
\mathrm{int}(\Theta)
%\reviseend
}\sum_{\omega\in\Omega_{n}}\|\nabla^{2}_{\theta}\log S_{\theta}(\omega)\|_{\mathrm{op}}\right),
\end{align*}
where $O(1)$ in the rightmost side implies a vector whose norm is $O(1)$.
This yields
\begin{align*}
&
\mathcal{G}_{\alpha}\left(\hat{\theta}^{(m)}_{\alpha}[\tilde{I}^{z,\omega^{*}}_{n}]\,;\,
\tilde{I}^{z,\omega^{*}}_{n}
\right)
-
\mathcal{G}_{\alpha}\left(\hat{\theta}^{(m)}_{\alpha}[\tilde{I}_{n}]\,;\,
\tilde{I}_{n}
\right)
\\
&=
\frac{\alpha}{(1-\alpha)n}
\sum_{\omega\in\Omega_{n}}
\Bigg{\{}
\frac{S_{\hat{\theta}^{(m)}_{\alpha}[\tilde{I}_{n}]}(\omega) }
{\alpha S_{\hat{\theta}^{(m)}_{\alpha}[\tilde{I}_{n}]}(\omega)
+
(1-\alpha)\tilde{I}_{n}(\omega)
}\\
&\qquad\qquad\qquad\qquad
-\frac{S_{\hat{\theta}^{(m)}_{\alpha}[\tilde{I}^{z,\omega^{*}}_{n}]}(\omega) }
{\alpha S_{\hat{\theta}^{(m)}_{\alpha}[\tilde{I}^{z,\omega^{*}}_{n}]}(\omega)
+
(1-\alpha)\tilde{I}^{z,\omega^{*}}_{n}(\omega)
}
\Bigg{\}}
\nabla_{\hat{\theta}^{(m)}_{\alpha}[\tilde{I}_{n}]}\log S_{\theta}(\omega)\\
&\quad
+O\left(\frac{2\varepsilon}{(1-\alpha)n}\sup_{\theta\in
%\mbox{\sout{$\Theta$}}
%\revisebegin
\mathrm{int}(\Theta)
%\reviseend
}\sum_{\omega\in\Omega_{n}}\|\nabla^{2}_{\theta}\log S_{\theta}(\omega)\|_{\mathrm{op}}\right)\\
%&=
%\frac{\alpha}{(1-\alpha)n}
%\sum_{\omega\in\Omega_{n}}
%\left\{
%\frac{S_{\hat{\theta}^{(m)}_{\alpha}[\tilde{I}_{n}]}(\omega) }
%{\alpha S_{\hat{\theta}^{(m)}_{\alpha}[\tilde{I}_{n}]}(\omega)
%+
%(1-\alpha)\tilde{I}_{n}(\omega)
%}
%-\frac{S_{\hat{\theta}^{(m)}_{\alpha}[\tilde{I}_{n}]}%(\omega) }
%{\alpha S_{\hat{\theta}^{(m)}_{\alpha}[\tilde{I}_{n}]}(\omega)
%+
%(1-\alpha)\tilde{I}^{z,\omega^{*}}_{n}(\omega)
%}
%\right.\\
%&\left.\qquad\qquad\quad\qquad\qquad
%+\frac{S_{\hat{\theta}^{(m)}_{\alpha}[\tilde{I}_{n}]}(\omega) }
%{\alpha S_{\hat{\theta}^{(m)}_{\alpha}[\tilde{I}_{n}]}(\omega)
%+
%(1-\alpha)\tilde{I}^{z,\omega^{*}}_{n}(\omega)
%}
%-\frac{S_{\hat{\theta}^{(m)}_{\alpha}[\tilde{I}^{z,\omega^{*}}_{n}]}(\omega) }
%{\alpha S_{\hat{\theta}^{(m)}_{\alpha}[\tilde{I}^{z,\omega^{*}}_{n}]}(\omega)
%+
%(1-\alpha)\tilde{I}^{z,\omega^{*}}_{n}(\omega)
%}
%\right\}
%\\
%&\qquad\qquad\quad\qquad\qquad
%\nabla_{\hat{\theta}^{(m)}_{\alpha}[\tilde{I}_{n}]}\log S_{\theta}(\omega)\\
%&\quad
%+O\left(\frac{2\varepsilon }{(1-\alpha)n}\sup_{\theta\in\Theta}\sum_{\omega\in\Omega_{n}}\|\nabla^{2}_{\theta}\log S_{\theta}(\omega)\|_{\mathrm{op}}\right)
%\\
&=\frac{\alpha}{(1-\alpha)n}(T_{1}+T_{2})+O\left(\frac{2\varepsilon }{(1-\alpha)n}
\sup_{\theta\in
%\mbox{\sout{$\Theta$}}
%\revisebegin
\mathrm{int}(\Theta)
%\reviseend
}\sum_{\omega\in\Omega_{n}}\|\nabla^{2}_{\theta}\log S_{\theta}(\omega)\|_{\mathrm{op}}\right),
\end{align*}
where 
\begin{align*}
T_{1}
&:=\sum_{\omega\in\Omega_{n}}
\Bigg{\{}
\frac{S_{\hat{\theta}^{(m)}_{\alpha}[\tilde{I}_{n}]}(\omega) }
{\alpha S_{\hat{\theta}^{(m)}_{\alpha}[\tilde{I}_{n}]}(\omega)
+
(1-\alpha)\tilde{I}_{n}(\omega)
}\\
&\qquad\qquad\quad-\frac{S_{\hat{\theta}^{(m)}_{\alpha}[\tilde{I}_{n}]}(\omega) }
{\alpha S_{\hat{\theta}^{(m)}_{\alpha}[\tilde{I}_{n}]}(\omega)
+
(1-\alpha)\tilde{I}^{z,\omega^{*}}_{n}(\omega)}
\Bigg{\}}\nabla_{\hat{\theta}^{(m)}_{\alpha}[\tilde{I}_{n}]}\log S_{\theta}(\omega),\\
T_{2}
&:=\sum_{\omega\in\Omega_{n}}
\Bigg{\{}\frac{S_{\hat{\theta}^{(m)}_{\alpha}[\tilde{I}_{n}]}(\omega) }
{\alpha S_{\hat{\theta}^{(m)}_{\alpha}[\tilde{I}_{n}]}(\omega)
+
(1-\alpha)\tilde{I}^{z,\omega^{*}}_{n}(\omega)
}\\
&\qquad\qquad\quad-\frac{S_{\hat{\theta}^{(m)}_{\alpha}[\tilde{I}^{z,\omega^{*}}_{n}]}(\omega) }
{\alpha S_{\hat{\theta}^{(m)}_{\alpha}[\tilde{I}^{z,\omega^{*}}_{n}]}(\omega)
+
(1-\alpha)\tilde{I}^{z,\omega^{*}}_{n}(\omega)
}\Bigg{\}} \nabla_{\hat{\theta}^{(m)}_{\alpha}[\tilde{I}_{n}]}\log S_{\theta}(\omega).
\end{align*}
Here we have
\begin{align*}
T_{1}&=\frac{\alpha S_{\hat{\theta}^{(m)}_{\alpha}[\tilde{I}_{n}]}(\omega^{*}) }{\alpha S_{\hat{\theta}^{(m)}_{\alpha}[\tilde{I}_{n}]}(\omega^{*})+
(1-\alpha)
\tilde{I}_{n}(\omega^{*})
}
\frac{(1-\alpha)z }{(1-\alpha) z + \alpha S_{\hat{\theta}^{(m)}_{\alpha}[\tilde{I}_{n}]}(\omega^{*})+
(1-\alpha)
\tilde{I}_{n}(\omega^{*})}\\
&\qquad\frac{\nabla_{\hat{\theta}^{(m)}_{\alpha}[\tilde{I}_{n}]}\log S_{\theta}(\omega^{*})}{\alpha}\\
&=
O\left(
\frac{1}{\alpha} \sup_{\theta,\,\omega}\| \nabla_{\theta}\log S_{\theta}(\omega) \|
\right)
\end{align*}
 and 
\begin{align*}
T_{2}&=\sum_{\omega\in\Omega_{n}}
\frac{(1-\alpha)\tilde{I}^{z,\omega^{*}}_{n}(\omega) }
{\alpha S_{\hat{\theta}^{(m)}_{\alpha}[\tilde{I}^{z,\omega^{*}}_{n}]}(\omega)+(1-\alpha)\tilde{I}^{z,\omega^{*}}_{n}(\omega) }\\
&\qquad\qquad
\frac{S_{\hat{\theta}^{(m)}_{\alpha}[\tilde{I}_{n}]}(\omega)-S_{\hat{\theta}^{(m)}_{\alpha}[\tilde{I}^{z,\omega^{*}}_{n}]}(\omega) }
{\alpha S_{\hat{\theta}^{(m)}_{\alpha}[\tilde{I}_{n}]}(\omega)+(1-\alpha)\tilde{I}^{z,\omega^{*}}_{n}(\omega)}
\nabla_{\hat{\theta}^{(m)}_{\alpha}[\tilde{I}_{n}]}\log S_{\theta}(\omega)\\
&=
O\left(
\frac{n\varepsilon}{\alpha}
\sup_{\theta,\theta',\omega}
\|\nabla_{\theta'}S_{\theta}(\omega)\|
\|(1/S_{\theta}(\omega))\nabla_{\theta}\log S_{\theta}(\omega)
\|
\right).
\end{align*}
Then we obtain
\begin{align*}
&
\left\|\mathcal{G}_{\alpha}\left(\hat{\theta}^{(m)}_{\alpha}[\tilde{I}^{z,\omega^{*}}_{n}]\,;\,
\tilde{I}^{z,\omega^{*}}_{n}
\right)
-
\mathcal{G}_{\alpha}\left(\hat{\theta}^{(m)}_{\alpha}[\tilde{I}_{n}]\,;\,
\tilde{I}_{n}
\right)
\right\|
\\
&=
\frac{1}{(1-\alpha)n} O(U_1)
+\varepsilon\Bigg\{
\frac{1}{(1-\alpha)} O(U_1^2)+
\frac{2}{(1-\alpha)} O(U_2) \Bigg\}
\end{align*}
and thus
\begin{align*}
&\left\|\hat{\theta}^{(m+1)}_{\alpha}[\tilde{I}_{n}^{z,\omega^{*}}]
-
\hat{\theta}^{(m+1)}_{\alpha}[\tilde{I}_{n}]\right\|\\
&\le
\frac{\gamma_{m}}{(1-\alpha)n} O(U_1)
+
\varepsilon\Bigg\{
1+
\frac{\gamma_{m}}{(1-\alpha)} O(U_1^2)
+ \frac{2\gamma_{m}}{(1-\alpha)} O(U_2) \Bigg\}.
\end{align*}
 This implies that the assertion holds for $m+1$ and by the mathematical induction, we get the conclusion.
\end{proof}

\section{Proof of Proposition \ref{prop: stability of line search}}
\label{appendix: Proof of line search}

%\proof
Let $\gamma$ be the step size satisfying 
the Armijo condition $\mathcal{A}[c\,,\,\theta\,,\,D_{\alpha}^{(n)}[\tilde{I}_{n} :S_{\theta}]]$
and $\gamma\le 2(1-c)/L$.
Let $\mathcal{P}_{\gamma}[\theta]$ be
\[
\mathcal{P}_{\gamma}[\theta]:=\theta-\gamma\nabla_{\theta}D_{\alpha}^{(n)}[\tilde{I}_{n}\,:\,S_{\theta}].
\]

We begin with utilizing the descent lemma (c.f., \cite{Bertsekas_NP}):
\begin{align}
D_{\alpha}^{(n)}
\left[\tilde{I}_{n}^{z,\omega^{*}}\,:\,
S_{\mathcal{P}_{\gamma}[\theta]}\right]
&\le 
D_{\alpha}^{(n)}
\left[\tilde{I}_{n}^{z,\omega^{*}}\,:\,
S_{\theta}\right]
+\langle
\nabla_{\theta}D_{\alpha}^{(n)}[\tilde{I}_{n}^{z,\omega^{*}}:S_{\theta}]
\,,\,
\mathcal{P}_{\gamma}[\theta]-\theta
\rangle\nonumber\\
&\quad +\frac{L}{2}\|\mathcal{P}_{\gamma}[\theta]-\theta\|^{2},
\label{eq: smoothness inequality}
\end{align}
where $\langle\cdot,\cdot \rangle$ denotes the Euclidean inner-product.
From equation (\ref{eq: comparison of first step}), we have
\begin{align*}
\left|\left\langle v\,,\, (\mathcal{P}_{\gamma}[\theta]-\theta) \right\rangle 
-\left\langle v\,,\, \left(-\gamma \nabla_{\theta}D_{\alpha}^{(n)}[\tilde{I}_{n}^{z,\omega^{*}}:S_{\theta}]\right)
\right\rangle\right|
&\le \frac{1}{n}\frac{\gamma}{1-\alpha}
\|v\|\|\nabla_{\theta}S_{\theta}(\omega^{*})\|\\
&\le \frac{1}{n}\frac{ 2U_{1}}{L(1-\alpha)}
\|v\|,
\end{align*}
where $U_{1}$ is defined in Assumption \ref{Assumption: norm constraint of spectral} and $L$ is defined in (\ref{eq: Lipschitz}). 
Putting the quantity $\nabla_{\theta}D_{\alpha}^{(n)}[\tilde{I}_{n}^{z,\omega^{*}}:S_{\theta}]$ into $v$ in the above inequality, we get
%This yields
\begin{align}
\langle
\nabla_{\theta}D_{\alpha}^{(n)}[\tilde{I}_{n}^{z,\omega^{*}}:S_{\theta}]
\,,\,
\mathcal{P}_{\gamma}[\theta]-\theta
\rangle
\le -\gamma
\|\nabla_{\theta}D_{\alpha}^{(n)}[\tilde{I}_{n}^{z,\omega^{*}}:S_{\theta}]\|^{2}
+\frac{1}{n}\frac{2U_{1}^{2}}{L(1-\alpha)^{2}},
\label{eq: approximate first}
\end{align}
where we use the inequality
\begin{align*}
    \|\nabla_{\theta}D_{\alpha}^{(n)}[\tilde{I}_{n}^{z,\omega^{*}}:S_{\theta}]\|
    \leq \frac{\alpha U_{1}}{(1-\alpha) n } \sum_{\omega\in\Omega_{n}}
    \left| \frac{(\alpha-1)S_{\theta}(\omega)+(1-\alpha)\tilde{I}_{n}^{z,\omega^{*}}(\omega)}{\alpha S_{\theta}(\omega)+(1-\alpha)\tilde{I}_{n}^{z,\omega^{*}}(\omega)}\right|
    \leq \frac{U_{1}}{1-\alpha}; 
\end{align*}
see \eqref{expr:grad-SpecRenyi} for the expression of $\nabla_{\theta}D_{\alpha}^{(n)}[\tilde{I}_{n}^{z,\omega^{*}}:S_{\theta}]$.
Equation (\ref{eq: comparison of first step}) also gives 
\begin{align*}
    \left\| \gamma\nabla_{\theta}D_{\alpha}^{(n)}[\tilde{I}_{n}^{z,\omega^{*}}:S_{\theta}] - \gamma\nabla_{\theta}D_{\alpha}^{(n)}[\tilde{I}_{n}:S_{\theta}] \right\|
    \leq \frac{\gamma U_{1}}{n(1-\alpha)}
    \leq \frac{2U_{1}}{nL(1-\alpha)}
\end{align*}
so that we obtain
\begin{align}
\|\mathcal{P}_{\gamma}[\theta]-\theta\|^{2}
\le 
\gamma^{2} \|\nabla_{\theta}D_{\alpha}^{(n)}[\tilde{I}_{n}^{z,\omega^{*}}:S_{\theta}]\|^{2}
+\left(\frac{2}{n}+\frac{1}{n^{2}}\right)\frac{4U_{1}^{2}}{L^{2}(1-\alpha)^{2}}.
\label{eq: approximate second}
\end{align}
Together with (\ref{eq: smoothness inequality}), equations (\ref{eq: approximate first}) and (\ref{eq: approximate second}) imply
\begin{align*}
D_{\alpha}^{(n)}
\left[\tilde{I}_{n}^{z,\omega^{*}}\,:\,
S_{\mathcal{P}_{\gamma}[\theta]}\right]
\le 
D_{\alpha}^{(n)}
\left[\tilde{I}_{n}^{z,\omega^{*}}\,:\,
S_{\theta}\right]
-\left(\gamma-\frac{L}{2}\gamma^{2}\right) \|\nabla_{\theta}D_{\alpha}^{(n)}[\tilde{I}_{n}^{z,\omega^{*}}:S_{\theta}]\|^{2}+\varepsilon_{n},
\end{align*}
where $\varepsilon_{n}$ is given by
\[
\varepsilon_{n}:=
\left(\frac{3}{n}+\frac{1}{n^{2}}\right)\frac{2U_{1}^{2}}{L(1-\alpha)^{2}}.
\]
For $\gamma$ satisfying $1-(L/2)\gamma \ge c$, we have
\begin{align*}
D_{\alpha}^{(n)}
\left[\tilde{I}_{n}^{z,\omega^{*}}\,:\,
S_{\mathcal{P}_{\gamma}[\theta]}\right]
&\le 
D_{\alpha}^{(n)}
\left[\tilde{I}_{n}^{z,\omega^{*}}\,:\,
S_{\theta}\right]
-c\gamma\|\nabla_{\theta}D_{\alpha}^{(n)}[\tilde{I}_{n}^{z,\omega^{*}}:S_{\theta}]\|^{2}+\varepsilon_{n}.
\end{align*}
Together with (\ref{eq: approximate second}), the assumption $\|\nabla_{\theta}D_{\alpha}^{(n)}[\tilde{I}_{n}:S_{\theta}]\|\ge \kappa$
implies that 
\[
\|\nabla_{\theta}D_{\alpha}^{(n)}[\tilde{I}_{n}^{z,\omega^{*}}:S_{\theta}]\|^{2}
\ge \kappa -\varepsilon_{n}/\underline{\gamma}^{2}.
\]
Thus, for sufficiently large $n$, we have $\kappa -\varepsilon_{n}/\underline{\gamma}^{2}>0$ so that
\begin{align*}
&D_{\alpha}^{(n)}
\left[\tilde{I}_{n}^{z,\omega^{*}}\,:\,
S_{\mathcal{P}_{\gamma}[\theta]}\right]\\
&\le 
D_{\alpha}^{(n)}
\left[\tilde{I}_{n}^{z,\omega^{*}}\,:\,
S_{\theta}\right]
-\left(c-\frac{\varepsilon_{n}}{\underline{\gamma}(\kappa-\varepsilon_{n}/\underline{\gamma}^{2})}\right)
\gamma\|\nabla_{\theta}D_{\alpha}^{(n)}[\tilde{I}_{n}^{z,\omega^{*}}:S_{\theta}]\|^{2},
\end{align*}
which concludes the proof.
\qed

\section{Proof of Proposition \ref{prop: failure Itakura Saito}}
\label{appendix: proof of Prop}

%\begin{proof}
The difference between $\hat{\theta}^{(1)}_{\mathrm{IS}}[\tilde{I}^{z,\omega^{*}}_{n}]$
and $\hat{\theta}^{(1)}_{\mathrm{IS}}[\tilde{I}_{n}]$
is explicitly written as
\begin{align*}
\|\hat{\theta}^{(1)}_{\mathrm{IS}}[\tilde{I}^{z,\omega^{*}}_{n}]-
\hat{\theta}^{(1)}_{\mathrm{IS}}[\tilde{I}_{n}]
\|
&=\|\{\hat{\theta}^{(1)}_{\mathrm{IS}}[\tilde{I}^{z,\omega^{*}}_{n}]-\theta^{(0)}\}
-
\{\hat{\theta}^{(1)}_{\mathrm{IS}}[\tilde{I}_{n}]
-\theta^{(0)}\}
\|\\
&=
\frac{\gamma_{1} z}{n}
\left\|\frac{\nabla_{\theta^{(0)}}\log S_{\theta}(\omega^{*}) }{S_{\theta^{(0)}}(\omega^{*})}
\right\|,
\end{align*}
which proves the first assertion.
For the second assertion,
observe
\begin{align*}
&\sum_{\omega\in\Omega_{n}}
\left[1-\frac{\tilde{I}^{z,\omega^{*}}_{n}(\omega)}{S_{\theta}(\omega)}\right]\nabla_{\theta}\log S_{\theta}(\omega)
\Bigg{|}_{\theta=\hat{\theta}^{(\infty)}_{\mathrm{IS}}[\tilde{I}^{z,\omega^{*}}_{n}] }
\\
&=
\sum_{\omega\in\Omega_{n}}
\left[1-\frac{\tilde{I}_{n}(\omega)}{S_{\theta}(\omega)}\right]\nabla_{\theta}\log S_{\theta}(\omega)
\Bigg{|}_{\theta=\hat{\theta}^{(\infty)}_{\mathrm{IS}}[\tilde{I}^{z,\omega^{*}}_{n}] }
-z\frac{\nabla_{\hat{\theta}^{(\infty)}_{\mathrm{IS}}[\tilde{I}^{z,\omega^{*}}_{n}]}
\log S_{\theta}(\omega^{*})
}
{S_{\hat{\theta}^{(\infty)}_{\mathrm{IS}}[\tilde{I}^{z,\omega^{*}}_{n}]}(\omega^{*}) }.
\end{align*}
This gives the second assertion, which completes the proof.
\qed

\section{Additional comparison of optimization paths}
\label{appendix: additional comparison}
%\revisebegin

This appendix provides an additional comparison of optimization paths.
Here we consider contaminated observations defined as
\[
X_{t}^{\circ}=X_{t}+\sqrt{z}\sin(t\pi/2), \quad t=1,\ldots,n=500
\]
with $(X_t)$ generated from AR(2) model having $(\sigma,\varphi_{1},\varphi_{2})$.

\begin{figure}[h]
    \centering
    \includegraphics[width=0.8\linewidth]{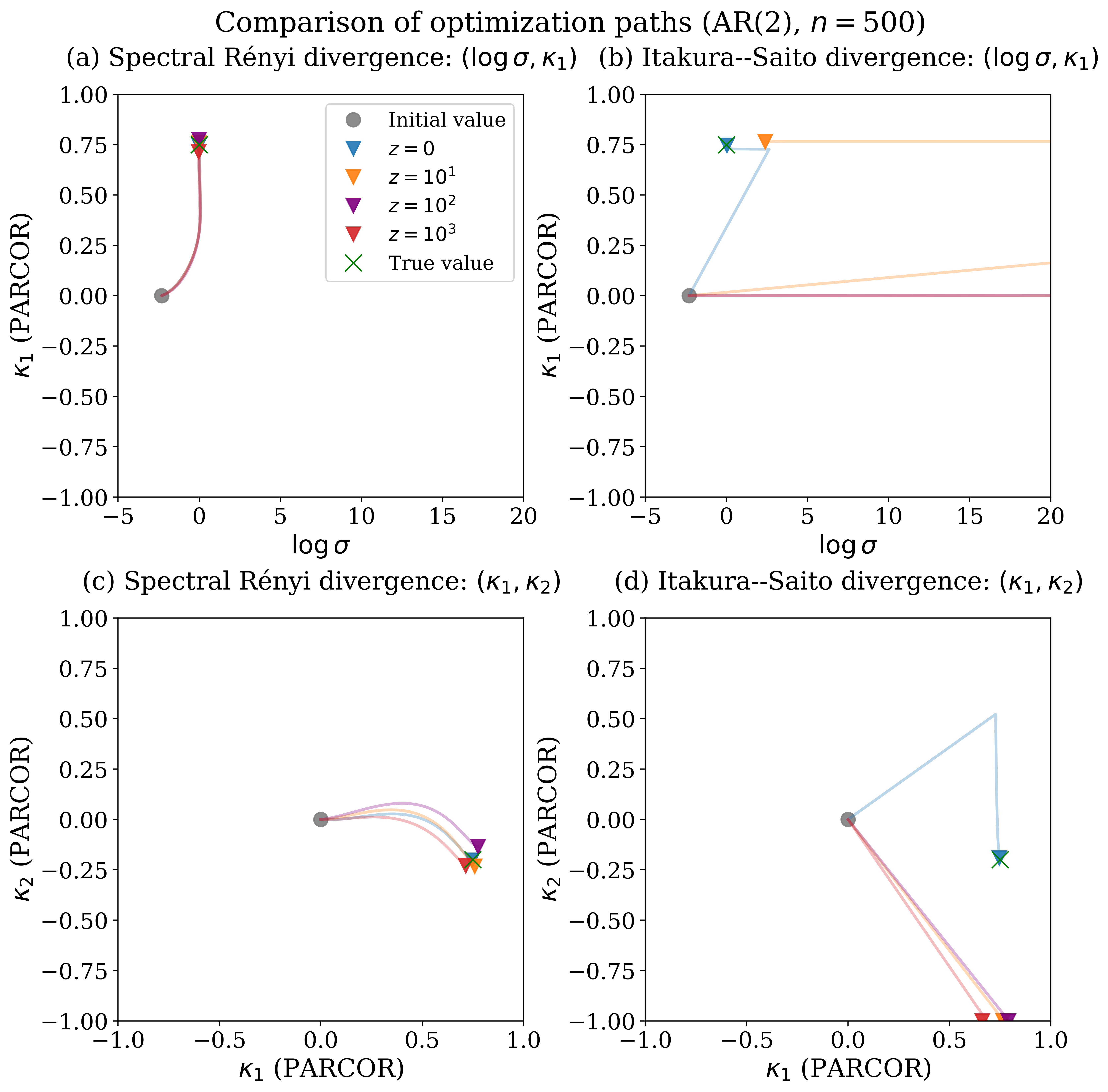}
    \caption{Comparison of optimization paths using the AR(2) model with the length $n$ of a time series set to $500$.
    Panels (a) and (c) show the optimization paths (curves with light colors) of the spectral R\'{e}nyi divergence minimization that start from the initial value (the gray circle) and terminate at points after $10000$ iterations (the inverted triangles), where the colors indicate different values of the contamination $z$.
    Panels (b) and (d) show the corresponding results for the Itakura--Saito divergence minimization. }
    \label{fig:trajectoryforAR(2)}
\end{figure}

To remove the parameter constraints arising from the stationarity condition, we employ the following transformation:
\begin{align}
\theta_{2}:= \,\mathrm{atanh}\left(\frac{\varphi_{1}}{1-\varphi_{2}}\right)\in\mathbb{R}
\quad\text{and}\quad
\theta_{3}:= \,\mathrm{atanh}(\varphi_{2}) \in \mathbb{R}. 
\label{eq: transformation}
\end{align}
This transformation is interpreted as follows.
First, we transform the AR coefficients $\varphi_{1},\varphi_{2}$ into the partial autocorrelation (PARCOR) coefficients $\kappa_{1},\kappa_{2}$ using the Durbin--Levinson formula:
\[
\kappa_{1}:=\frac{\varphi_{1}}{1-\varphi_{2}}\in (-1,1)
\quad\text{and}\quad
\kappa_{2}:=\varphi_{2} \in (-1,1).
\]
Next, we map the PARCOR coefficients to unrestricted real values via the Fisher transformation:
\[
\theta_{2}:= \,\mathrm{atanh}\left(\kappa_{1}\right)\in\mathbb{R}
\quad\text{and}\quad
\theta_{3}:= \,\mathrm{atanh}(\kappa_{2}) \in \mathbb{R},
\]
which yields the transformation in (\ref{eq: transformation}).
We perform 10000 iterations of gradient descent with a step size of 0.01 with respect to $(\log \sigma, \theta_{2},\theta_{3})$.

Figure \ref{fig:trajectoryforAR(2)} shows a comparison of the optimization paths for the AR(2) model, where the true parameter values are set to $(\sigma^{*},\varphi^{*}_{1},\varphi^{*}_{2})=(1,0.9,-0.2)$
and the hyperparameter $\alpha$
is set to $\alpha = 0.9$.
Even for the PARCOR coefficients (in addition to $\sigma$), the optimization trajectory based on the Itakura--Saito divergence remains sensitive to outliers, whereas that based on the spectral R\'enyi divergence remains stable with respect to all parameters.

%\reviseend

\section{Additional simulation studies}
\label{appendix: additional examples}

This appendix presents additional numerical experiments.

Tables \ref{tab: Brune without trend additional}
and \ref{tab: Brune with trend additional} show the estimation results for the Brune spectral model with attenuation in Section \ref{subsec: Brune}
on the basis of the different initial values 
$\theta^{(0),4}=(0.1,0.1,0.1)$
and $\theta^{(0),5}=(2,2,2)$.
For almost all the settings of initial values, the spectral density based on the spectral $\alpha=0.5$-R\'{e}nyi divergence performs the best.

\begin{table}[ht]
\caption{
The mean values of biases with standard deviations without any trend. For each initial value, the values closest to zero are underlined. R\'{e}nyi is abbreviated as R; Itakura--Saito is abbreviated as IS; Initial value is abbreviated as Init. Values greater than 10 are denoted by $*$.
}
\centering
\begin{tabular}{r|crrr}
  & Init & $\hat{\sigma}-\sigma^{*}$ & $\hat{\omega_{c}}-\omega_{c}^{*}$ & $\hat{Q}-Q^{*}$ \\ 
  \hline
R($\alpha=0.50$) & 
$\theta^{(0),4}$ &\underline{-0.03} ($\pm$0.05) & \underline{-0.11} ($\pm$0.07) & \underline{0.26} ($\pm$0.15) \\ 
  R ($\alpha=0.75$) & $\theta^{(0),4}$ & 
0.04 ($\pm$0.06) & -0.30 ($\pm$0.03) &  1.33 ($\pm$0.17)  \\ 
  R ($\alpha=0.90$) &$\theta^{(0),4}$& 
0.09 ($\pm$0.08) & -0.39 ($\pm$0.03) &  6.13 ($\pm$0.04) \\ 
  IS with $I_{n}$  & $\theta^{(0),4}$ & 
  $*$
  %$5\cdot 10^{15}$ ($\pm 8\cdot 10^{14}$) 
  & 
  $*$
  %$1\cdot 10^{16}$ ($\pm 1\cdot 10^{15}$) 
  &  
  $*$
  %$8\cdot 10^{16}$ ($\pm 1\cdot 10^{16}$) 
  \\ 
  IS with $I_{n}^{S}$ & $\theta^{(0),4}$ & 
  $*$
%$5\cdot 10^{18}$ ($\pm9\cdot 10^{18}$)  
& 
$*$
%$1\cdot 10^{19}$ ($\pm 1\cdot 10^{19}$) 
&
$*$
%$9\cdot 10^{19}$ ($\pm 1\cdot 10^{20}$) 
\\ 
   \hline  \hline
R ($\alpha=0.50$) & 
$\theta^{(0),5}$ & \underline{-0.01}($\pm$0.06) & -0.18 ($\pm$0.05) & 0.47 ($\pm$0.16)\\ 
  R ($\alpha=0.75$) & $\theta^{(0),5}$ & \underline{0.01} ($\pm$ 0.06) & -0.21 ($\pm$0.05) & 0.68 ($\pm$0.22) \\ 
  R ($\alpha=0.90$) &$\theta^{(0),5}$& 0.02 ($\pm$0.06) & -0.24 ($\pm$0.04) & 1.18 ($\pm$0.34)\\ 
  IS with $I_{n}$  & $\theta^{(0),5}$ & 0.07 ($\pm$0.11) & \underline{-0.08} ($\pm$0.12) & \underline{0.12} ($\pm$0.18)\\ 
  IS with $I_{n}^{S}$ & $\theta^{(0),5}$ & -0.28 ($\pm$0.16) & 1.84 ($\pm$1.66) & 1.90 ($\pm$0.63)\\ 
  \hline
\end{tabular}
\label{tab: Brune without trend additional}
\end{table}

\begin{table}[ht]
\caption{
The mean values of biases with standard deviations with the trigonometric trends. For each initial value, the values closest to zero are underlined. R\'{e}nyi is abbreviated as R; Itakura--Saito is abbreviated as IS; Initial value is abbreviated as Init.
{%\revisebegin
Values greater than 10 are denoted by $*$.
%\reviseend
}
}
\centering
\begin{tabular}{r|crrr}
  & Init & $\hat{\sigma}-\sigma^{*}$ & $\hat{\omega_{c}}-\omega_{c}^{*}$ & $\hat{Q}-Q^{*}$ \\ 
  \hline
R ($\alpha=0.50$) & 
$\theta^{(0),4}$ & \underline{0.21} ($\pm$0.07) & \underline{-0.21} ($\pm$0.06) & \underline{0.24}  ($\pm$0.14)\\ 
  R ($\alpha=0.75$) & $\theta^{(0),4}$ & 
0.61 ($\pm$0.09) & -0.43 ($\pm$0.02) &  1.14 ($\pm$0.20)  \\ 
  R ($\alpha=0.90$) &$\theta^{(0),4}$& 
1.65 ($\pm$0.12) & -0.60 ($\pm$0.02) &  6.05 ($\pm$0.05) \\ 
  IS with $I_{n}$  & $\theta^{(0),4}$ 
  &
  $*$
  %$5\cdot 10^{15}$ ($\pm 8\cdot 10^{14}$) 
  &
  $*$
  %$1\cdot 10^{16}$ ($\pm 1 \cdot 10^{15}$) 
  &
  $*$
  %$8\cdot 10^{16}$ ($\pm 1\cdot 10^{16}$) 
  \\ 
  IS with $I_{n}^{S}$ & $\theta^{(0),4}$ 
  &
  $*$
%$5\cdot 10^{18}$ ($\pm 9 \cdot 10^{18}$) 
&
$*$
%$1\cdot 10^{19}$ ($\pm 1 \cdot 10^{19}$) 
&
$*$
%$9\cdot 10^{19}$ ($\pm 1\cdot 10^{20}$)
\\ 
   \hline  \hline
R ($\alpha=0.50$) & 
$\theta^{(0),5}$ & \underline{0.23} ($\pm$0.07) &  \underline{-0.26} ($\pm$0.04) & 0.41 ($\pm$0.16)\\ 
  R ($\alpha=0.75$) & $\theta^{(0),5}$ & 0.55 ($\pm$0.08) &  -0.35 ($\pm$0.04) & 0.55 ($\pm$0.20) \\ 
  R ($\alpha=0.90$) &$\theta^{(0),5}$ & 1.50 ($\pm$0.10) &  -0.50 ($\pm$0.03) & 0.85 ($\pm$0.35)\\ 
  IS with $I_{n}$  & $\theta^{(0),5}$ & 3.91 ($\pm$0.07) &  -0.34 ($\pm$0.02) & \underline{-0.40} ($\pm$0.02)\\ 
  IS with $I_{n}^{S}$ & $\theta^{(0),5}$ & 2.85 ($\pm$0.95) &  -0.17 ($\pm$0.51) & 2.26 ($\pm$0.93)\\ 
  \hline
\end{tabular}
\label{tab: Brune with trend additional}
\end{table}

\begin{figure}[ht]
    \centering
\includegraphics[scale=0.45]{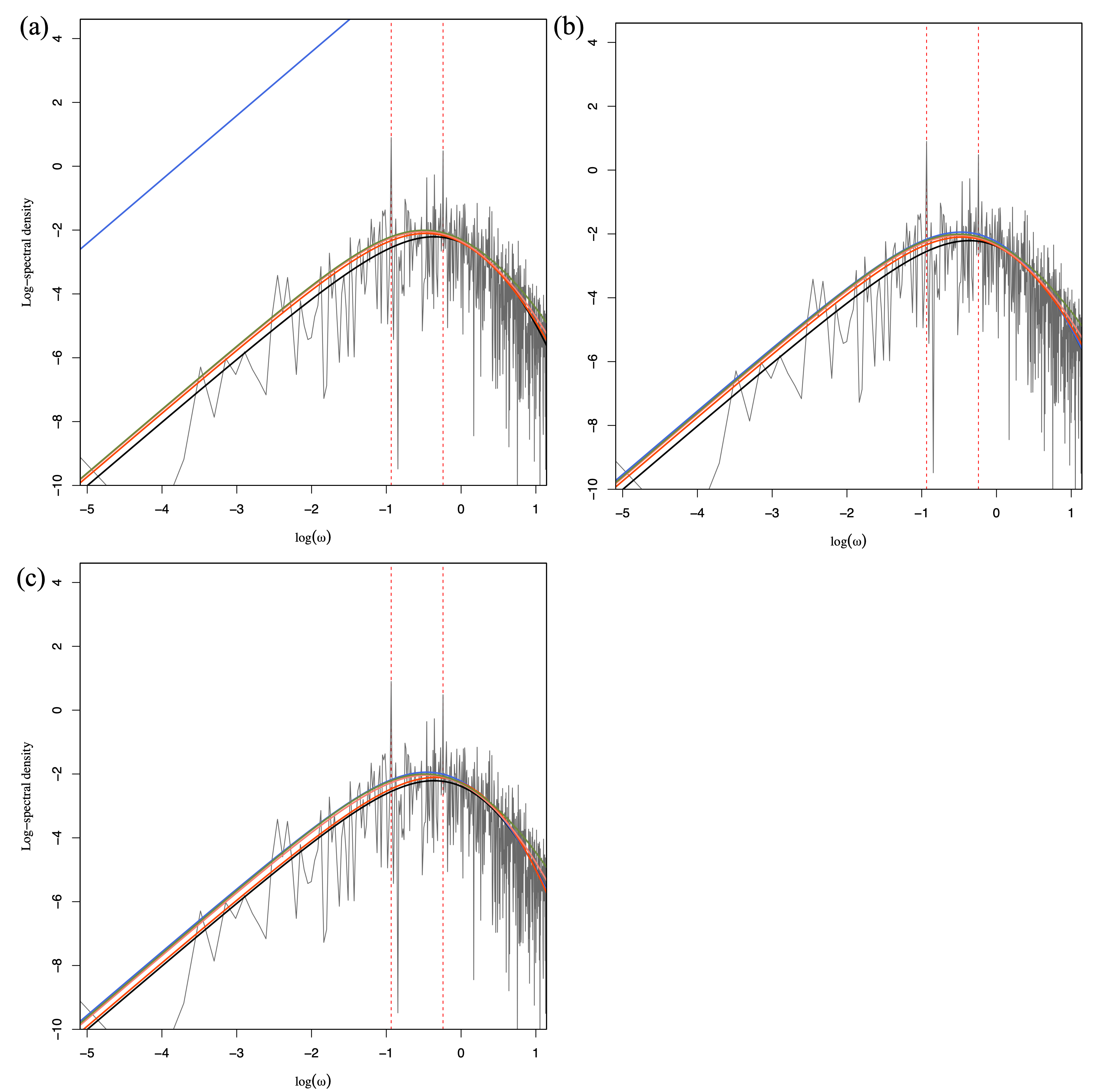}
    \caption{Spectral densities with the estimates plugged in for the Brune spectral model with attenuation and with the trigonometric trends of $z_{1}=z_{2}=2.5$. The gray curve is the periodogram. The true spectral density is colored in black. Spectral densities based on the spectral R\'{e}nyi divergence ($\alpha=0.5,0.75,0.9$) are colored in red, salmon pink, and green,respectively. The spectral density based on the Itakura--Saito divergence is colored in blue.
    The red dashed lines denote the frequencies at which the outliers are injected.
    (a) the result based on the initial value $\theta^{(0),1}$,
    (b) the result based on the initial value $\theta^{(0),2}$,
    (c) the result based on the initial value $\theta^{(0),3}$.
    }
    \label{Fig_SPD_Brune_z10}
\end{figure}

\begin{figure}[ht]
    \centering
\includegraphics[scale=0.45]{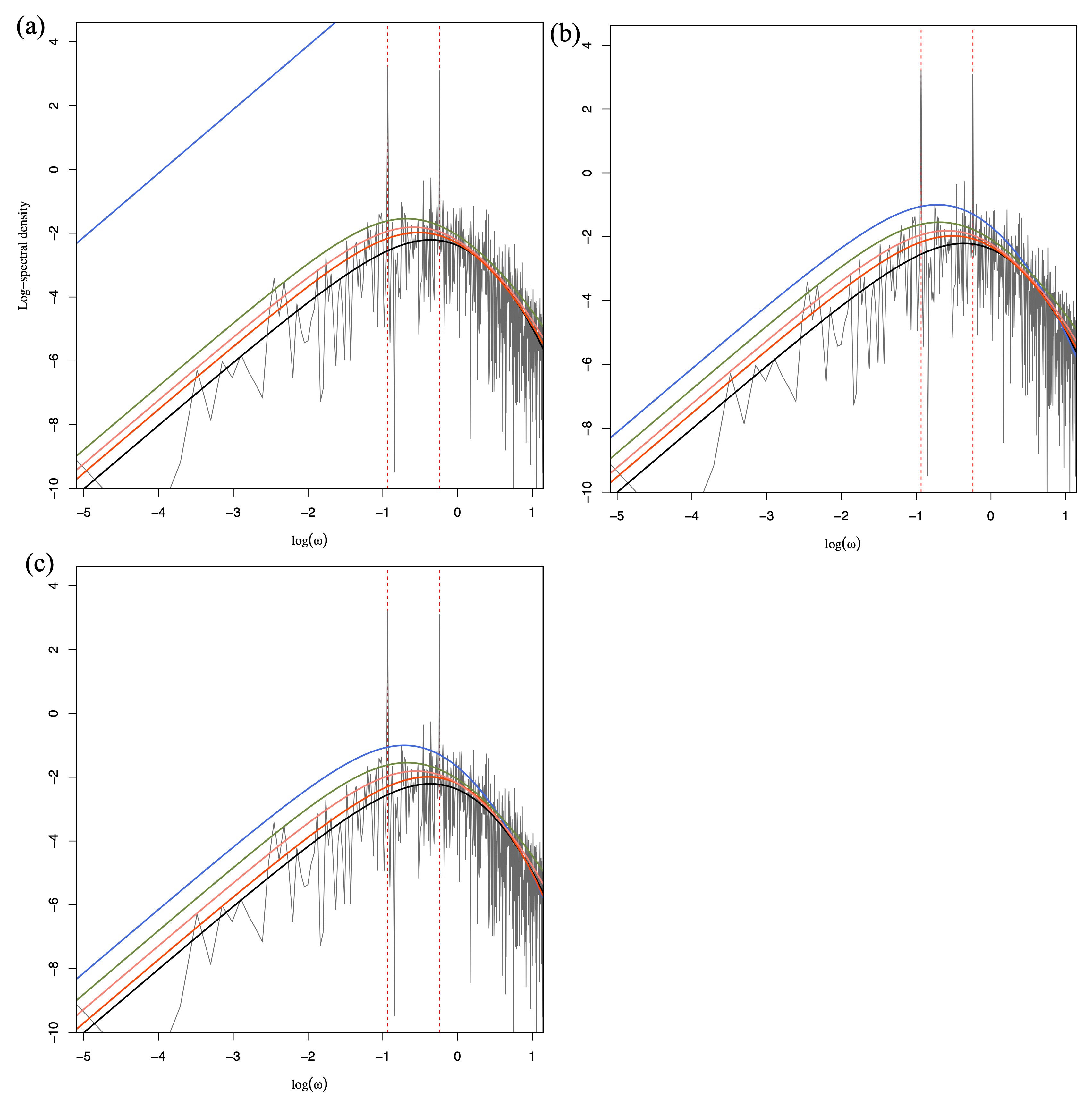}
    \caption{Spectral densities with the estimates plugged in for the Brune spectral model with attenuation and with the trigonometric trends of $z_{1}=z_{2}=25$. The gray curve is the periodogram. The true spectral density is colored in black. Spectral densities based on the spectral R\'{e}nyi divergence ($\alpha=0.5,0.75,0.9$) are colored in red, salmon pink, and green,respectively. The spectral density based on the Itakura--Saito divergence is colored in blue. The red dashed lines denote the frequencies at which the outliers are injected. (a) the result based on the initial value $\theta^{(0),1}$,
    (b) the result based on the initial value $\theta^{(0),2}$,
    (c) the result based on the initial value $\theta^{(0),3}$.
    }
    \label{Fig_SPD_Brune_z100}
\end{figure}

Figures \ref{Fig_SPD_Brune_z10}
and \ref{Fig_SPD_Brune_z100} display an instance of a set of the estimated spectral densities.
In each figure, the following representations are used:
\begin{itemize}
\item The solid gray curve represents the periodogram $I_{n}$;
\item the solid black curve denotes the true spectral density $S^{\mathrm{BA}}_{\theta^{*}}$;
\item the solid blue curve illustrates the spectral density with the minimum Itakura--Saito divergence estimate plugged in;
\item the solid red, salmon pink, and green curves display the spectral density with the minimum spectral R\'{e}nyi divergence ($\alpha=0.9,0.75,0.5$) estimates  plugged in, respectively.
\end{itemize}

For any initial value and for any strength of trends, the spectral density based on the spectral $\alpha=0.5$-R\'{e}nyi divergence performs the best.

\clearpage

\bibliographystyle{spbasic_mod} 
\bibliography{specdiv}

@article{liu1989limited,
  title={On the limited memory method for large scale optimization},
  author={D. Liu and J. Nocedal},
  journal={Mathematical Programming},
  series={B},
  volume={45},
  number={3},
  pages={503--528},
  year={1989},
  doi={10.1007/BF01589116}
}

@inproceedings{malitsky2020adaptive,
  title={Adaptive Gradient Descent Without Descent},
  author={Y. Malitsky and K. Mishchenko},
  booktitle={Proceedings of the 37th International Conference on Machine Learning},
  volume={119},
  pages={6702--6712},
  year={2020},
  publisher={PMLR}
}

@misc{yagishita2025simple,
  title={Simple linesearch-free first-order methods for nonconvex optimization},
  author={S. Yagishita and M. Ito},
  howpublished={\url{https://arxiv.org/abs/2509.14670}},
  note={arXiv:2509.14670},
  year={2025}
}

@article{kano2025spatiotemporal,
  title = {Spatio-temporal characteristics in the {GEONET} {F5} solution in the frequency domain estimated based on the robust spectral analysis},
  author = {M. Kano and K. Yano and Y. Tanaka and T. Takabatake and Y. Ohta},
  journal = {Earth, Planets and Space},
  volume = {77},
  number = {103},
  year = {2025},
  doi = {10.1186/s40623-025-02236-3}
}

@article{AttouchBolteSvaiter2013,
author={H. Attouch and J. Bolte and B. Svaiter},
title={Convergence of descent methods for semi-algebraic and tame problems: proximal algorithms, forward–backward splitting, and regularized {G}auss–{S}eidel methods},
journal={Mathematical Programming},
year={2013},
volume={137},
pages={91--129}
}

@article{Kurdyka_1998,
author={K. Kurdyka},
title={On gradients of functions definable in $o$-minimal structures},
journal={Annales de l’institut Fourier},
volume={48},
year={1998},
pages={769--783}
}

@incollection{Lojasiewicz_1963,
author={S. {\L}ojasiewicz},
year={1963},
title = {Une propri{\'e}t{\'e} topologique des sous-ensembles analytiques r{\'e}els},
  booktitle = {Les {\'E}quations aux {D}{\'e}riv{\'e}es Partielles},
  publisher = {{\'E}ditions du Centre National de la Recherche Scientifique},
  address = {Paris},
  pages = {87--89}
}

@book{Nocedal_Wright_book,
  title={Numerical optimization},
  author={Nocedal, Jorge and Wright, Stephen J},
  year={2006},
  edition={Second},
  publisher={Springer},
  address={New York}
}

@book{Maronna2019robust,
  title={Robust statistics: theory and methods (with {R})},
  author={Maronna, Ricardo A and Martin, R Douglas and Yohai, Victor J and Salibi{\'a}n-Barrera, Mat{\'\i}as},
  year={2019},
  edition={second},
  publisher={John Wiley \& Sons},
  address={Hoboken.}
}

@book{Diakonikolas_Kane_2023,
  title={Algorithmic high-dimensional robust statistics},
  author={Diakonikolas, Ilias and Kane, Daniel M},
  year={2023},
  publisher={Cambridge University Press},
  address={Cambridge}
}

@book{Bertsekas_NP,
  title={Nonlinear Programming},
  author={Bertsekas, D. P.},
  year={1999},
  edition={second},
  publisher={Athena Scientific},
  address={Belmont, Massachusetts}
}

@article{Kakizawaetal1998,
author = {Y. Kakizawa and R. Shumway and M. Taniguchi},
title = {Discrimination and Clustering for Multivariate Time Series},
journal = {Journal of the American Statistical Association},
volume = {93},
number = {441},
pages = {328--340},
year = {1998},
publisher = {Taylor \& Francis},
doi = {10.1080/01621459.1998.10474114}
}

@article{Hirukawa2005,
  title={Cluster analysis for non-{G}aussian locally stationary processes},
  author={Hirukawa, Junichi},
  journal={International Journal of Theoretical and Applied Finance},
  volume={9},
  number={01},
  pages={113--132},
  year={2006},
  publisher={World Scientific}
}

@article{Langbein_2004,
author={J. Langbein},
year={2004},
title={Noise in two-color electronic distance meter measurements revisited},
journal={Journal of Geophysical Research: Solid Earth},
volume={109}
}

@article{Kolmogorov_1941,
author={A. Kolmogorov},
title={Stationary sequences in {H}ilbert space},
journal={Bulletin of Moscow State University, Mathematics},
year={1941},
pages={1--40}
}

@article{Szego_1920,
author={O. Szeg\"{o}},
title={Beitr\"{a}ge zur Theorie der {T}oeplitzschen Formen},
journal={Mathematische Zeitschrift},
volume={6},
year={1920},
pages={167--202} 
}

@article{ZhangTaniguchi1995,
  title={Nonparametric approach for discriminant analysis in time series},
  author={Zhang, Guoqiang and Taniguchi, Masanobu},
  journal={Journal of Nonparametric Statistics},
  volume={5},
  number={1},
  pages={91--101},
  year={1995},
  publisher={Taylor \& Francis}
}

@incollection{Martin2005,
  title={Statistical methods for the enhancement of noisy speech},
  author={Martin, Rainer},
  booktitle={Speech Enhancement},
  chapter={2},
  pages={43--65},
  year={2005},
  publisher={Springer},
  address={Berlin}
}

@article{Brune1970,
author={J. Brune},
title={Tectonic stress and spectra of seismic shear waves from earthquakes},
journal={Journal of Geophysical Research},
volume={75},
year={1970},
pages={4997--5009}
}

@book{Matern1960,
author={B. Mat\'{e}rn},
year={1960},
title={Spatial Variation: Stochastic Models and Their Application to Some Problems in Forest Surveys
and Other Sampling Investigations},
publisher={Stockholm: Statens Skogsforskningsinstitut}
}

@misc{Shayevitz2010,
author={P. Shayevitz},
title={A note on a characterization of {R}\'{e}nyi measures and its relation to composite hypothesis testing},
year={2010},
note={arXiv:1012.4401v1}
}

@article{vanErvenHarremos,
author={T. van Erven and P. Harremo{\"e}s},
title={R{\'e}nyi divergence and {K}ullback--{L}eibler divergence},
journal={IEEE Transactions on Information Theory},
volume={60},
pages={3797--3820},
year={2014}
}

@book{AkiandRichard,
author={K. Aki and P. G. Richards},
year={1980},
title={Quantitative seismology: theory and methods},
publisher={W.~H.~Freeman and Co.},
address={San Francisco}
}

@article{Giletal2013,
author={M. Gil and F. Alajaji and T. Linder},
journal={Information Sciences},
volume={249},
title={{R}\'{e}nyi divergence measures for commonly used univariate continuous distributions},
year={2013},
pages={124--131}
}

@article{HeydeandDai1996,
author={C. Heyde and W. Dai},
title={ON THE ROBUSTNESS TO SMALL TRENDS OF ESTIMATION BASED ON THE SMOOTHED PERIODOGRAM},
journal={Journal of Time Series Analysis},
volume={17},
year={1996},
pages={141--150}
}

@article{McCloskeyandPerron2013,
year={2013},
author={A. McCloskey and P. Perron},
title={MEMORY PARAMETER ESTIMATION IN THE PRESENCE OF LEVEL SHIFTS AND DETERMINISTIC TRENDS},
journal={Econometric Theory},
volume={29},
pages={1196--1237}
}

@article{Iacone2010,
author={F. Iacone},
year={2010},
journal={{\it Journal of Time Series Analysis}},
volume={31},
pages={37-49}, 
title={Local {W}hittle estimation of the memory parameter in presence of deterministic components}
}

@article{Griveletal2021,
    author = {E. Grivel and R. Diversi and F. Merchan},
    title = {{K}ullback--{L}eibler and {R}\'{e}nyi divergence rate for {G}aussian stationary {ARMA} processes comparison},
    journal = {Digital Signal Processing},
    year = {2021},
    volume={116},
    pages={103089}
}

@incollection{Parzen1993,
  author    = {Parzen, E.},
  title     = {Stationary Time Series Analysis Using Information and Spectral Analysis},
  booktitle = {Developments in Time Series Analysis. In Honour of
M. B. Priestley},
  pages     = {139--148},
  year      = {1993},
  publisher = {Chapman \& Hall},
  address   = {London},
  editor={S. Rao},
}

@article{Whittle1953,
title={The Analysis of Multiple Stationary Time Series},
author={P. Whittle},
journal={Journal of the Royal Statistical Society. Series B (Methodological)},
volume={15}, 
year={1953},
pages={125--139}
}

@book{Vajda1989,
title={Theory of Statistical Inference and Information},
author={I. Vajda},
year={1989},
publisher={Springer},
address={Dordrecht}
}

@article{FujisawaEguchi2008,
author={H. Fujisawa and S. Eguchi},
year={2008},
title={Robust parameter estimation with a small bias against heavy contamination},
journal={Journal of Multivariate Analysis},
volume={99},
pages={2053--2081}
}

@article{CalderoniandAbercrombie2023,
  title={Investigating spectral estimates of stress drop for small to moderate earthquakes with heterogeneous slip distribution: {E}xamples from the 2016--2017 {A}matrice earthquake sequence},
  author={Calderoni, Giovanna and Abercrombie, Rachel E},
  journal={Journal of Geophysical Research: Solid Earth},
  volume={128},
  number={6},
  pages={e2022JB025022},
  year={2023}
}

@article{Yoshimitsuetal2023,
title={Estimation of source parameters using a non-{G}aussian probability density function in a {B}ayesian framework},
author={N. Yoshimitsu and T. Maeda and T. Sei},
journal={Earth, Planets and Space},
volume={75},
pages={33},
year={2023}
}

@book{BrockwellandDavis,
author={P. Brockwell and R. Davis},
title={Time Series: Theory and Methods},
edition={2nd},
publisher={Springer},
year={1991},
address = {New York}
}

@article{Taniguchi1987,
author={M. Taniguchi},
title={Minimum Contrast Estimation for Spectral Densities of Stationary Processes},
journal={Journal of the Royal Statistical Society. Series B },
year={1987},
pages={315--325},
volume={49}
}

@article{ItakuraSaito1968,
author={F. Itakura and S. Saito},
year={1968},
title={Analysis synthesis telephony based on the maximum likelihood method},
journal={Proceedings of the 6th of the International Congress on Acoustics},
pages={C17--C20}
}

@book{Pinsker1964,
author={M. Pinsker},
title={Information and information stability of random variables and processes},
note={Translated and annotated by A. Feinstein},
year={1964},
publisher={Holden-Day Inc},
address={San Francisco}
}

\end{document}